\documentclass[11pt,a4paper, oneside]{article}
\usepackage[utf8]{inputenc}
\usepackage[english]{babel}
\usepackage{amsfonts}
\usepackage{amsmath}
\usepackage{amssymb}
\usepackage{amsthm}
\usepackage{comment}
\usepackage{float}
\usepackage{textcomp}
\usepackage{faktor}
\usepackage{setspace}
\usepackage[dvipsnames]{xcolor}
\usepackage{graphicx}
\usepackage{fancybox}
\usepackage{mathrsfs}
\usepackage{mathtools}
\usepackage{tikz}
\usepackage{centernot}
\usepackage{esint}
\usepackage{algpseudocode}
\usepackage{framed}
\usepackage{hyperref}
\usepackage[bottom]{footmisc}
\usepackage[square, sort, numbers]{natbib}
\usepackage[T1]{fontenc}
\usepackage{mathtools}
\usepackage{graphicx}
\usepackage{epstopdf}
\usepackage{listings}
\usepackage{fancyhdr}
\usepackage{afterpage}
\usepackage{ dsfont }
\usepackage{soul}
\usepackage[a4paper,top=3.5cm,bottom=3.5cm,left=3cm,right=3cm]{geometry}
\tikzstyle{level 1}=[level distance=4cm, sibling distance=3.5cm,->]
\tikzstyle{level 2}=[level distance=4cm, sibling distance=2cm,->]

\tikzstyle{bag} = [text width=2em, text centered]
\tikzstyle{end} = []

\DeclarePairedDelimiter{\set}{ \{ } { \} }

\newcommand{\sfrac}[2]{ 
	{\raise0.8ex\hbox{$#1$} \!\mathord{\left/ 
			{\vphantom {#1 #2}}\right.\kern-\nulldelimiterspace} 
		\!\lower0.8ex\hbox{$#2$}} 
}

\newtheorem{teo}{Theorem}[section]
\newtheorem{lema}[teo]{Lemma}
\newtheorem{co}[teo]{Corollary}
\newtheorem{prop}[teo]{Proposition}

\theoremstyle{definition}
\newtheorem{defin}[teo]{Definition}
\newtheorem{notz}{Notation}
\newtheorem{obs}[teo]{Remark}

\newcommand{\R}{\mathbb{R} }
\newcommand{\E}{\mathbb{E} }
\newcommand{\bbG}{\mathbb{G} }
\newcommand{\bbX}{\mathbb{X} }
\newcommand{\bbF}{\mathbb{F} }
\newcommand{\bbU}{\mathbb{U}}
\newcommand{\F}{\mathcal{F}}
\newcommand{\calG}{\mathcal{G}}
\newcommand{\calA}{\mathcal{A}}
\newcommand{\calU}{\mathcal{U}}
\newcommand{\calL}{\mathcal{L}}
\newcommand{\calE}{\mathcal{E}}
\newcommand{\caltildeE}{\widetilde{\mathcal{E}}}
\newcommand{\calH}{\mathcal{H}}
\newcommand{\calR}{\mathcal{R}}
\newcommand{\calD}{\mathscr{D}}
\newcommand{\hatu}{\hat{u}}
\newcommand{\hatX}{\hat{X}}
\newcommand{\hatY}{\hat{Y}}

\newcommand{\hatq}{\hat{q}}
\newcommand{\hatk}{\hat{\kappa}}
\newcommand{\hatb}{\hat{b}}

\newcommand{\hatg}{\hat{g}}

\newcommand{\hatp}{\hat{p}}

\newcommand{\tH}{\widetilde{H}}

\newcommand{\probspace}{(\Omega,\F,P)}
\renewcommand{\t}{t \in [0,T]}

\DeclareMathSymbol{\shortminus}{\mathbin}{AMSa}{"39}

\date{20th February 2023}
\author{Giulia di Nunno\thanks{Department of Mathematics, University of Oslo, P:O: Box 1053 Blindern, N-0316 Oslo, Email:giulian@math.uio.no.} \thanks{Department of Business and Management Science, NHH Norwegian School of Economics, Helleveien 30, N-5045 Bergen.}\and Michele Giordano\thanks{Department of Mathematics, University of Oslo, P:O: Box 1053 Blindern, N-0316 Oslo, Email:michelgi@math.uio.no}}

\title{Maximum principles for \\ stochastic time-changed Volterra games}

\begin{document}
\thispagestyle{empty}
\numberwithin{equation}{section}
\allowdisplaybreaks
\maketitle
\begin{abstract}
	\noindent	We study a stochastic differential game between two players, controlling a forward stochastic Volterra integral equation (FSVIE). Each player has to optimize his own performance functional which includes a backward stochastic differential equation (BSDE). The dynamics considered are driven by time-changed Lévy noises, with absolutely continuous time-change process. We prove a sufficient maximum principle to characterize Nash equilibria and the related optimal strategies. For this we use techniques of control under partial information, and the non-anticipating stochastic derivative. The zero-sum game is presented as a particular case.

\vspace{1mm}

\end{abstract}

\section{Introduction}
Hereafter we consider a type of stochastic games with memory. The goal of each one of the two players $i = 1,2$ is to find an optimal control $u_i$ that maximizes his personal performance functional
\begin{equation}\label{J_intro}
	J_i(u):=\E\left[\int_0^TF_i(t,u(t),X(t),Y_i(t)) dt+\varphi_i(X(T))+\psi_i(Y_i(0))\right],
\end{equation} 
which includes profit rate $F_i$, utility function $\phi_i$ and $\psi_i$, associated to recursive utility or risk evaluation. Here $u=(u_1,u_2)$, with $u_i$ representing the control of player $i$. The controlled process $X$ follows the forward dynamics of Volterra type :
\begin{equation}\label{X_t intro}
X(t)=X_0+\int_0^t b(t,s,\lambda_s,u(s),X(s\shortminus))ds+\int_0^t\!\!\!\int_\R \kappa(t,s,z,\lambda_s,u(s),X(s\shortminus))\mu(dsdz),  \ \t,
\end{equation} 
while, for each $i=1,2$, $Y_i$ is a backward stochastic differential equation (BSDE) with dynamics given by:
\begin{align}
	Y_i(t)&=h_i(X(T))-\int_t^Tg_i(s,\lambda_s,u(s),X(s\shortminus),Y_i(s\shortminus),\Theta_i(s,\cdot))ds\nonumber \\
	&\quad +\int_t^T\!\!\!\int_\R\Theta_i(s,z)\mu(dsdz)+\int_t^TdM_i(s), \quad \t \label{Y_t intro}.
\end{align}
The noise $\mu$, appearing in \eqref{X_t intro} and \eqref{Y_t intro}, is the mixture of a conditional Brownian motion and a conditional centred Poisson random measure. They can be regarded as a time-changed Lévy noises with rate $\lambda$. The inclusion of time-change allows to gain more flexibility in modelling than what is possible with the classical Lévy framework, as pointed out, e.g., in \cite{Swishchuk,BNielsen}. Indeed, time change has been suggested in the study of volatility modelling \cite{CW, CGMY, DiNunno1,Yablonski}, energy models \cite{CW}, default models \cite{OksendalGames} and also kinetic theory \cite{Tusheng-Oksendal-Sulem}. In this work, we deal with a time-change process of the form 
\begin{equation}\label{timechange}
	\Lambda_t=\int_0^t\lambda_sds, \quad \t,
\end{equation}
where $\lambda$ is the stochastic time-change rate. In this case, the noise $\mu$ may loose the independent increment structure proper of Lévy noise and the Markovianity, thus allowing for time dependencies. The above term $M_i$, ($i=1,2$) is a martingale orthogonal to $\mu$ naturally appearing in the martingale representation theorem when working with noises and information flows that do not have the \emph{perfect} stochastic integral representation property, see e.g. \cite{Shiryaev}. Indeed, we recall that Kunita-Watanabe result shows that, given a square integrable martingale $\mathcal{M}$, as integrator, and its natural filtration
\begin{equation*}
    \bbF^{\mathcal{M}} = \left\{\F^{\mathcal{M}}_t, \ \t\right\},
\end{equation*}
any square integrable $\F_T^{\mathcal{M}}$-measurable random variable $\xi$ admits representation
\begin{equation}\label{perfect integral repr}
    \xi = \xi^{\perp} + \int \varphi \  d\mathcal{M}
\end{equation}
by means of a unique stochastic integrand $\varphi$. Here $\xi^{\perp}$ is a stochastic remainder orthogonal to all the stochastic integrals with respect to $\mathcal{M}$. It is well known that $\xi^{\perp}$ is a constant (naturally equal to $\E[\xi]$) whenever $\mathcal{M}$ is a Gaussian or a centered Poisson random measure, or a mixture of the two, and the reference filtration is generated by $\mathcal{M}$, see e.g. \cite{Int_stoc, chou-meyer, DiNunnoRF, DiNunnoStoch}. We say that the representation \eqref{perfect integral repr} is \emph{perfect} if $\xi^{\perp}$ is a constant. In this case, the space of all stochastic integrals with respect to $\mathcal{M}$ coincides with $L^2(\Omega, \F_T^{\mathcal{M}}, P)$ modulo constants. This is tightly connected with chaos expansions (see e.g. \cite{DOP-book, DiNunnoRF}). The study of general BSDEs of type \eqref{Y_t intro} can be found in e.g. \cite{Papantoleon}.

Each BSDEs in \eqref{Y_t intro} corresponds to one of the players in the stochastic game and is connected to the FSVIE \eqref{X_t intro} in the sense that they depend on $X$ and are driven by the same kind of noise. 

In view of the dynamics \eqref{X_t intro}, this work deals with stochastic games with memory. Some examples of papers that consider games with memory are \cite{Shi-Zhu, FBSDEgames, Shi-Wang2}. In both \cite{FBSDEgames} and \cite{Shi-Zhu} though, the memory effect is obtained by means of a delay whereas in \cite{Shi-Wang2}, similarly to what we do here, the authors consider a forward equation of Volterra type, but work with linear-quadratic stochastic integral games.

In terms of structure of the game studied, our work is in line with the stochastic differential games in as \cite{OksendalGames, FBGames, FBSDEgames}, where the authors consider games based on one forward stochastic differential equation and two backward stochastic differential equations. 
Other examples of stochastic games in the literature are the \textit{purely forward} stochastic games (see, e.g., \cite{Carmona,Wang-Yu,TTK} to cite a few) where there are no backward equations such as \eqref{Y_t intro} and the \textit{purely backward} stochastic games (see \cite{Shi-Zhu}) where there is no forward equation such as \eqref{X_t intro}. By considering a system of forward-backward stochastic differential equations we are able to have a setting that is comprehensive of both cases above. Notice that our setting is also more general than the one in the above mentioned \cite{OksendalGames, FBGames, FBSDEgames} as we consider the forward equation \eqref{X_t intro} to be of Volterra type and deal with backward equations with an additional martingale term which is due to the more general nature of the noise considered. Lastly, we remark that the structure of the game presented here is different from e.g. \cite{Shi-Wang2, Chen-Yong}, where the authors consider systems of forward-backward stochastic differential equations, but exploit the backward equations as a tool to express conditions on the existence of Nash equilibria for the game.

The first element of novelty of this paper with respect to, e.g., \cite{OksendalGames,FBSDEgames,FBGames, Shi-Zhu}, is the introduction of the time-changed Lévy process $\mu$ as a driver for both the forward and backward dynamics. This allows us to go beyond the classical Brownian motion and pure Lévy framework, by considering some quite general, but still treatable martingales. Moreover, contrarily to \cite{OksendalGames, FBGames,Shi-Zhu} we consider forward dynamics with memory, introduced by means of Volterra coefficients, similarly to what is presented in \cite{Oksendal1, Shi-Wang1} and used there for describing coupled systems.

Another novelty introduced in the current work comes from considering the BSDEs under an enlargement of the filtration generated by $\mu$. Then we shall see how the different filtrations introduce different integral representations. This will be exploited in the definition of the adjoint backward equations associated to the maximum principle approach.


Last but not least, differently from \cite{FBSDEgames} for instance, we do not make use of Malliavin derivatives. In fact, dealing with time-changed processes, we cannot directly apply the standard Malliavin calculus nor its Hida-Malliavin extension, see e.g. \cite{DOP-book}, because the driving noise is not Brownian nor of Poisson type. Still one may think of applying a conditional form of the Malliavin calculus, as presented in \cite{Yablonski}, but this would require to impose some additional condition on the domain of our processes, conditions that would be difficult to characterize since they would depend on the control process itself. At the current time there is no Hida-Malliavin type extension for conditional Malliavin calculus available, hence we cannot use either this approach in our work as it was instead done in the framework of \cite{Oksendal1}. For this reason we make use of the  \textit{non-anticipating derivative} (NA-derivative) with Itô-type calculus. The NA-derivative, introduced in \cite{DiNunnoDerivatives} for general martingales and then extended to martingale random fields in \cite{DiNunno2}, is directly connected to the Itô integral and provides explicit stochastic integral representations. Moreover, being the domain of the NA-derivative the whole $L^2(dP)$, we are able to overcome the domain restriction issues that would emerge from using the (conditional) Malliavin calculus. See e.g. \cite{DiNunno2,DiNunnoDerivatives} and also \cite{DiNunnoGiordano} where the NA-derivative is used in the context of optimal portfolios.

 We also put notice that stochastic Volterra games are connected with forward-backward stochastic Volterra integral equations (FBSVIEs). Even though some authors provided some interesting results concerning FBSVIEs (see e.g. \cite{Oksendal1,Chen-Yong,Shi-Wang1}), the results in this paper stand out from those in bibliography by providing a different setting. This is achieved both by considering equations driven by a time-changed Lévy noise and by working under different information flows.

In our work, we use partial information techniques. In fact, we work under a filtration $\bbF\subseteq \bbG$ regarded as \emph{partial} with respect to a filtration $\bbG$, containing in addition the future values of the time-change process $\Lambda$. Also, similarly to e.g. \cite{OksendalGames,TTK, FBGames}, we assume that the individual information flows available to each players are different. This is represented by taking two sub-filtrations $\calE^{(i)}_t\subseteq \F_t$, $i=1,2$. 

We provide sufficient maximum principle for non-zero sum stochastic games and we present a zero-sum game reformulation of these results as a particular case.

This paper is organized as follows: in Section 2 we present the framework for the study of time-changed Lévy process and we recall the concept of NA-derivative. In Section 3 we describe the optimization problem and provide a sufficient maximum principle for non-zero sum games under both the information flow $\bbF$ and $\bbG$. As a particular case we also provide the zero-sum game reformulation of such theorem. In Section 4 we provide some utility consumption examples both for the zero and non zero-sum games and show that in some simple cases we can actually manage to derive the optimal control for our forward-backward game. Lastly, in Section 5, we make some conclusive remarks on the use of Volterra-type BSDEs in place of \eqref{Y_t intro}.

\section{Time change and the framework}
We adopt the framework introduced in \cite{DiNunno1} and also exploited in \cite{DiNunnoGiordano}. Hereafter we summarize the crucial points fundamental for this study. 
Let us consider a complete probability space $\probspace$ and define $$\bbX:=[0,T]\times \R	=\Big([0,T]\times \set{0}\Big)\cup\Big([0,T]\times \R_0\Big),$$ where $\R_0:=\R\backslash \set{0}$ and $\mathcal{B}_\bbX$ is the Borel $\sigma$-algebra on $\bbX$. We define $\mathcal{L}$, the set of stochastic processes $\lambda=(\lambda^B,\lambda^H)$ such that each component $l=B,H$ is a non-negative, stochastically continuous process in $L^1(P\times dt)$.
Let $\nu$ be a $\sigma$-finite measure on the Borel sets $\mathcal{B}_{\R_0}$ of $\R_0$ satisfying $\int_{\R_0}z^2\nu(dz)<\infty$.
We define the random measure $\Lambda$ on $\mathcal{B}_\bbX$ by
\begin{equation*}
\Lambda(\Delta):=\int_0^T\mathds{1}_{\set{(t,0)\in\Delta}}(t)\lambda_t^B dt+\int_0^T \!\!\!\int_{\R_0}\mathds{1}_\Delta(t,z)\nu(dz)\lambda_t^H dt, \quad \Delta \subseteq \bbX,
\end{equation*}
and denote the $\sigma$-algebra generated by all the values of $\Lambda$ by $\F^\Lambda$. Also, we define the filtration $\bbF^\Lambda:=\set{\F_t^\Lambda, \t}$ with $F_t^\Lambda$ generated by the values $\Lambda(\Delta\cap [0,t])$, $\Delta\subseteq\bbX$. Set $\F^\Lambda:=\F_T^\Lambda$.

\begin{defin}
	The signed random measure $\mu$ on $\mathcal{B}_\bbX$ is defined as the mixture 
	\begin{equation}\label{definitionmu}
	\mu(\Delta):=B\left(\Delta\cap[0,T]\times\set{0}\right)+\tH\left(\Delta\cap[0,T]\times\R_0\right), \quad \Delta \subseteq \bbX,
	\end{equation}
	of $B$ and $\tH$, which are a conditional Gaussian measure (given $\F^\Lambda$) and a doubly stochastic centered Poisson random measure (given $\F^\Lambda$), respectively. The random measures $B$ and $\tH$ are indeed connected to time-changed L\'evy noises as shown in, e.g., \cite{Serfozo}, \cite[Theorem 3.1]{Grigelionis}.
\end{defin}

In the sequel we consider two types of information flows. The control problem is naturally associated with the filtration generated by the values $\mu(\Delta)$,  $\Delta\subset[0,t]\times \R$, $\t$, i.e. $\bbF^\mu:=\set{\F_t^\mu, \ \t}$. 
Since the filtration $\bbF^\mu$ might not be right-continuous, we consider the information flow 
$$\bbF:=\set{\F_t, \ \t}, \quad \text{ where } \F_t:=\bigcap_{r>t}\F_r^\mu.$$ 
The second information flow of interest is 
$$\bbG:=\set{\calG_t, \ \t}, \quad \text{ where }  \calG_t:=\F_t^\mu\vee\F^\Lambda.$$ 
Notice that $\bbG$ is right-continuous and that, while at the horizon $T$ we have $\calG_T=\F_T$, the initial information is different in the two flows: $\calG_0=\F^\Lambda$ and $\F_0$ is trivial. Namely $\bbG$ is enlarged of the information on the future values of $\Lambda$. In the sequel we technically exploit the interplay between the two filtrations regarding $\bbF$ as partial with respect to $\bbG$.

The random measure $\mu$ is a martingale random field with values in $L^2(dP)$ with respect to both filtration $\bbG$ and $\bbF$, see \cite{DiNunno2} and \cite[Remark 2.6]{DiNunno1}. Given a predictable random field $\phi$ in $L^2(d\Lambda\times dP)$, we have that
\begin{equation}\label{QuadVar}
	\E\bigg[\left(\iint_\bbX \phi(s,z) \mu(dsdz)\right)^2\bigg]=\E\bigg[\iint_\bbX \phi(s,z)^2\Lambda(dsdz)\bigg],
\end{equation}
under both the information flows $\bbG$ and $\bbF$.
We denote $\mathcal{I}^{\bbG}$ (resp. $\mathcal{I}^{\bbF}$) the subspace of $L^2(d\Lambda\times dP )$ of the random fields admitting a $\bbG$-predictable (resp. $\bbF$-predictable) modification.

	When considering the (It\^o type) stochastic integration we have the following integral representation with an explicitly characterized integrand by means of the non-anticipating (NA) derivative. This result first appeared in \cite{DiNunnoDerivatives} (for general $L^2$-valued martingales as integrators) and its extension to the random fields in \cite{DiNunno2}. In case of integration with respect to the information flow $\bbG$, the statement reads as follows.
	\begin{teo}\label{C-Ok-nonanticipating}
	For any $\xi\in L^2(dP)$ the NA-derivative $\calD \xi$ is well defined in $\mathcal{I}^\bbG$ as the $L^2(d\Lambda \times dP)$ limit  
$\calD \xi = \lim_{n \to \infty} \varphi_n$ 	
	of the simple predictable random fields $(\varphi_n)_{n\in \mathbb{N}}$ in the form
	$$
	\varphi_n(t,z) := \sum_{k=1}^{K_n} \E \bigg[ \xi \frac{\mu(\Delta_{nk})}{ \E [\Lambda(\Delta_{nk}) \vert \calG_{s_{nk}}]} \Big\vert \calG_{s_{nk}} \bigg] \mathbf{1}_{\Delta_{nk}}(t,z), \quad (t,z) \in \bbX. 
	$$
	defined on some intersecting system $(\Delta_{nk})_{k=1,...,K_n; n\in \mathbb{N}} $ with $\Delta_{nk} = (s_{nk}, u_{nk}] \times B_{nk}$ and $B_{nk}$ belonging to a countable semiring generating $\mathcal{B}_\R$. Furthermore, the following stochastic integral representation of $\xi$ holds
	\begin{equation}\label{C-O_non_ant}
		\xi=\xi^{0}+\iint_\bbX\calD_{t,z} \xi  \ \mu(dtdz),
	\end{equation}
	where the stochastic term $\xi^{0}\in L^2(dP)$ is such that  $\calD\xi^{0}\equiv 0$ and in particular $\xi^{0}=\E[\xi|\F^\Lambda]$. 
\end{teo}
\noindent The duality between the It\^o stochastic integral and the NA-derivative is an immediate consequence.
	\begin{co}\label{Yab} For all $\phi$ in $\mathcal{I}^{\bbG}$ and all $\xi$ in $L^2(dP)$, we have
	\begin{equation}\label{duality}
		\E\left[\xi\iint_\bbX\phi(t,z)\mu(dtdz)\right]=\E\left[\iint_\bbX \phi(t,z)\calD_{t,z}\xi \ \Lambda(dtdz)\right].
	\end{equation}
\end{co}
\noindent We observe that representation \eqref{C-O_non_ant} and the duality \eqref{duality} would hold also for stochastic integration with respect to the filtration $\bbF$. In that case, however, we would still have $\calD\xi^{0}=0$, but \emph{no explicit} representation for $\xi^{0}$.

\section{Maximum principles for time-changed Volterra games}
Having in mind a financial context, the FSVIE \eqref{X_t intro} represents the market value, whereas each BSDE \eqref{Y_t intro} is associated to some risk evaluation of the market as perceived by player $i=1,2$ or to a recursive utility in the stochastic game we consider.
For this reason the BSDEs will depend on the values of $X$, with dynamics
\begin{equation}\label{X_t}
	X(t)=X_0+\int_0^t b(t,s,\lambda_s,u(s),X(s\shortminus))ds+\int_0^t\!\!\!\int_\R \kappa(t,s,z,\lambda_s,u(s),X(s\shortminus))\mu(dsdz),
\end{equation}
where $X(0)=X_0\in\R$ and
\begin{align*}
	b&:[0,T]\times[0,T]\times[0,\infty)^2\times\bbU\times\R\times\Omega\longrightarrow\R,\\
	\kappa&:[0,T]\times[0,T]\times\R\times[0,\infty)^2\times\bbU\times\R\times\Omega\longrightarrow\R,
\end{align*}
are given functions with $\bbU$ specified in Definition \ref{defUoptim} below. 
We assume that $b$ and $\kappa$ are $C^2$ with respect to their first variable and, for all $t \in[0,T]$, $\lambda \in [0,\infty)^2$, $u\in\bbU$, $x\in\R$, the stochastic processes $s\longmapsto b(t,s,\lambda,u,x)$ and $s\longmapsto\kappa(t,s,z,\lambda,u,x)$, with $s \leq t$, are $\bbF$-predictable. 
We also assume that $b$ and $\kappa$ are $C^2$ with respect to $u$ and $x$, with partial derivatives in $L^2(dt\times dP)$ and $L^2(d\Lambda\times dP)$, respectively. 
 Later on we consider the coefficients $b$ and $\kappa$ in a functional setup:
\begin{align*}
	b&:[0,T]\times[0,T]\times\Xi_{\R^2_+}\times \Xi_{\bbU}\times\Xi_\R\times \Omega\longrightarrow\R, \\
	\kappa&:[0,T]\times[0,T]\times \R \times\Xi_{\R^2_+}\times \Xi_{\bbU}\times\Xi_\R\times \Omega\longrightarrow\R,
	\end{align*}
 where we denoted by $\Xi_S$ the space of measurable function on $[0,T]$ with values in $S$. Then we can interpret the coefficients in \eqref{X_t} via the evaluation at the point $s\in[0,T]$:
 \begin{align*}
     b(t,\cdot,\lambda_\cdot, u(\cdot), X^u(\cdot))(s) &= b(t,s,\lambda_s,u(s),X^u(s\shortminus))\\
    \kappa(t,\cdot,z,\lambda_\cdot, u(\cdot), X^u(\cdot))(s) &= \kappa(t,s,z,\lambda_s,u(s),X^u(s\shortminus)).
 \end{align*}
 We assume that $b$ and $\kappa$ are Fréchet differentiable (in the standard topology of càdlàg paths) with $C^2$ regularity in $t$, $x$ and $u$ (with the corresponding derivatives). In general, throughout the paper, $\partial_x$ denotes the appropriate partial derivative with respect to the variable $x$.

Lastly, we assume that, for all $z\in\R$, $\lambda \in [0,\infty)^2$, $u\in\calU$, $x\in\R$ the partial derivative of $\kappa$ with respect to $t$ (denoted with $\partial_t \kappa(t,s,z,\lambda,u,x)$) is locally bounded (uniformly in $t$) and satisfies
	\begin{equation}\label{HPkap}
	    |\partial_t\kappa(t_1,s,z,\lambda,u,x)- \partial_t \kappa(t_2,s,z,\lambda,u,x)|\leq K |t_1-t_2|,
	\end{equation}
	for some $K>0$ and for each fixed $s\leq t$, $\lambda \in [0,\infty)^2$, $u\in\calU$, $x\in\R$.
In this case, in fact, we are able to apply the Transformation Rule (see \cite{Protter} for the original result and see \cite{DiNunnoGiordano} Lemma 3.4 for its formulation in the context of time-change) and write $X$ in differential notation as
\begin{align}
    dX(t) &= \Bigg[ b(t,t,\lambda_t ,u(t), X(t)) + \int_0^t\partial_t b(t,s,\lambda_s, u(s), X(s)) ds \label{dX_t}\\
    &\quad + \int_0^t\!\!\!\int_\R \partial_t \kappa(t,s,z,\lambda_s, u(s) X(s) \mu(dsdz)\Bigg] dt + \int_\R\kappa(t,t,\lambda_t,u(t), X(t))\mu(dtdz)\nonumber.
\end{align}

The actions of the players $i=1,2$ are conveyed into the controls $u_i$, which form the stochastic process $u(s)=(u_1(s),u_2(s))$, $s \in [0,T]$. 

\noindent Note that if $b$ and $\kappa$ are Lipschitz continuous in $x$ and have at most linear growth in $x$, uniformly in $t,s\in[0,T]^2$, $\lambda\in[0,\infty)^2$, $ u\in \bbU$, then  $X$ admits an $\bbF$-adapted solution in $L^2(dt\times dP)$. For details concerning the existence and uniqueness of the solution of \eqref{X_t} we refer to Theorem 4.2 in \cite{DiNunnoGiordano}.

\begin{obs}\label{remark funzionale}
    We notice that by rewriting \eqref{X_t} in the differential form \eqref{dX_t}, we can link our Volterra optimization problem to a \emph{functional} SDE optimization problem. 
    As stated also in \cite{DiNunnoGiordano} Remark 3.3 some existence results for \eqref{dX_t} are available (see e.g. \cite{OksendalDhal, DiNunnoBanos}) and could possibly be extended to the current framework. In that case some additional conditions on $b$ and $\kappa$ are expected since none of those deals with time-changed Lévy noises $\mu$ as driver.
\end{obs}

For each player $i$, we consider the controlled BSDEs $(Y_i,\Theta_i ,M_i)$, depending on \eqref{X_t} of the form:
\begin{align}
	Y_i(t)&=h_i(X(T))-\int_t^Tg_i(s,\lambda_s,u(s),X(s\shortminus),Y_i(s\shortminus),\Theta_i(s,\cdot))ds\nonumber\\
	&\quad +\int_t^T\!\!\!\int_\R\Theta_i(s,z)\mu(dsdz)+\int_t^T dM_i(s), \quad \t,\label{Y_t}
\end{align}
where $M_i\in \mathcal M^{\perp}$, which is the space of the $\bbF$-martingales in $L^2(dt\times dP)$ orthogonal to $\mu$. 
For all $i=1,2$, 
\begin{align*}
h_i&:\R\longrightarrow\R\\
	g_i&:[0,T]\times[0,\infty)^2\times\bbU\times\R^2\times\mathcal{Z}\times\Omega\longrightarrow\R,
\end{align*}
where $\mathcal{Z}$ is the space of functions $\theta:\R\longrightarrow\R$ such that
\begin{equation*}
	|\theta(0)|^2+\int_{\R_0}\theta(z)^2\nu(dz)<\infty.
\end{equation*}

We assume that, for $i=1,2$,  the processes $s\longmapsto g_i(t,\lambda,u,x,y,\theta)$ is $\bbF$-adapted, integrable in $L^2(dt\times dP)$ for all $\lambda\in [0,\infty)^2$, $u\in \mathbb{U}$, $x\in \R$, $y\in\R$, $\theta \in \mathcal{Z}$, $C^2$ with respect to $t$, $x$ and $u$, and with partial derivatives in $L^2(dP)$. Also we assume that the functions $h_i:\R\longrightarrow\R $ are $\F_T$-measurable and $C^1$ with respect to $x$.

Throughout the paper we assume that for all $i=1,2$, there exists a solution (under $\bbF$) to \eqref{Y_t} $(Y_i,\Theta_i,M_i)\in L^2(dt\times dP)\times \mathcal{I}^{\bbF}\times \mathcal M^{\perp}$.

As for the coefficients of $X$ in \eqref{X_t}, we will also consider $g_i$, $i=1,2$ in a functional setting:
\begin{equation*}
    g_i:[0,T]\times \Xi_{\R^2}\times \Xi_{\bbU}\times \Xi_{\R}\times \calR^{\bbG} \times \mathcal Z\times\Omega\longrightarrow\R,
\end{equation*}
where $\mathcal{R}^\bbG$ is the space of all $\bbG$-adapted stochastic processes in $L^2(dt\times dP)$. We also require $g_i$, $i=1,2$ is Fréchet differentiable (in the standard topology of càdlàg paths) with $C^2$ regularity in $t$ $x$ and $u$ (with the corresponding derivatives).
Sufficient conditions that ensure the existence of a solution of \eqref{Y_t} in $L^2(dt\times dP)\times \mathcal{I}^{\bbF}\times \mathcal M^{\perp}$ with respect to a general filtration have been studied in \cite{Papantoleon}.

\vspace{5 mm}

Each player $i$ intends to optimize his position according to his own performance functional consisting of a utility evaluation and a risk assessment or recursive utility.
Indeed, let $F_i:[0,T]\times \bbU\times\R^2\longrightarrow\R$, $\varphi_i:\R\longrightarrow\R$ and $\psi_i:\R\longrightarrow\R$ be three functions that can be regarded as profit rates, utility function and recursive utility functions or risk evaluations of player $i$, respectively. We define the performance functional of each player $i=1,2$ as
	\begin{equation}\label{J}
		J_i(u):=\E\left[\int_0^TF_i(t,u(t),X(t),Y_i(t)) dt+\varphi_i(X(T))+\psi_i(Y_i(0))\right],
	\end{equation} 
	provided the integrals and expectations exist. We assume that $\varphi_i$ and $\psi_i$ are $C^1$ with respect to $x$, and that $F_i(\cdot,u,x,y)$ is $\bbF$-adapted for all $u,x,y$, and $C^1$ with respect to $u,x,y$, for all $\t$.
In this study, we allow for players to access different levels of information and we thus assume that the filtration of player $i$ is given by the right-continuous sub-filtrations $\E^{(i)}=\{\calE^{(i)}_t, \ \t\}$, with $\calE^{(i)}_t\subseteq \F_t$.	
\begin{defin}\label{defUoptim}
	For each player $i=1,2$, the set $\calA^{\calE}_i$ of admissible controls is given by  the $\E^{(i)}$-predictable processes with values in $A_i\subseteq \R^d$, $d\geq 1$ such that \eqref{X_t} and \eqref{Y_t} have a solution,  $h_i(X(T)) \in L^2(dP)$ for all $i=1,2$ and
	\begin{align*}
	    \E&\Bigg[\int_0^T|F_i(t,u(t),X(t),Y_i(t))| dt+|\varphi_i(X(T))|+|\psi_i(Y_i(0))| +\\
	    &|\partial_x\varphi_i(X(T))|^2 + |\partial_y\psi_i(Y_i(0))|^2 \Bigg]<\infty.
	\end{align*} We define $\bbU:=A_1\times A_2$.
\end{defin}

\begin{defin}\label{Nash}
	A Nash equilibrium for the forward-backward stochastic game \eqref{X_t}-\eqref{Y_t}-\eqref{J} is a pair $(\hatu_1,\hatu_2)\in\calA^{\calE}_1\times\calA^{\calE}_2$ such that
	\begin{align*}
		J_1(u_1,\hatu_2)&\leq J_1(\hatu_1,\hatu_2) \text{ for all } u_1\in\calA^{\calE}_1,\\
	J_2(\hatu_1,u_2)&\leq J_2(\hatu_1,\hatu_2) \text{ for all } u_2\in\calA^{\calE}_2.
	\end{align*}
	\end{defin}
\noindent	Intuitively this means that at equilibrium one player has no incentive in changing his strategy as long as the other player does not deviate from its own.

	\begin{obs}\label{EnlargedE}
	In the context of anticipated information, one could define a Nash equilibrium with respect to the enlarged filtration $\bbG$.
In particular, when working under $\bbG$, we consider an initial enlargement also for the information available to each player, i.e. $\tilde\E^{(i)}:=\E^{(i)}\vee \bbF^\Lambda =\{\caltildeE^{(i)}_t:=\calE^{(i)}_t\vee \F^\Lambda, \ \t\}$. 
Correspondingly, we define the set of admissible control processes for player $i$ for the stochastic game with anticipating information as $\calA_i^{\caltildeE}$ in line with Definition \ref{defUoptim}. Hence the Nash equilibrium of Definition \ref{Nash} will be defined for $(\hatu_1,\hatu_2)\in\calA^{\caltildeE}_1\times\calA^{\caltildeE}_2$.
We summarize the relationship among information flows in the chart below:
	$$\begin{matrix}
		&\F_t & \subset & \calG_t&:=&\F_t\vee \F^\Lambda\\
		&\rotatebox[origin=c]{90}{$\subset$} & &\rotatebox[origin=c]{90}{$\subset$} & &\\
		&\calE_t & \subset &\caltildeE_t&:=&\calE_t\vee \F^\Lambda.
	\end{matrix}$$
\end{obs}

\vspace{2mm}
In the sequel we work under the assumption that a Nash equilibrium for the forward-backward stochastic game \eqref{X_t}-\eqref{Y_t} with optimization functionals \eqref{J} exists. Our goal is to present results that allow us to find such equilibrium. In view of the non-Markovian nature of these systems and the general nature of the driving noise, we consider the maximum principle approach.

We start by defining the Hamiltonians functionals $\calH_i$ associated with the optimization problem of each player: 
\begin{equation*}
	\calH_i:[0,T]\times\Xi_{\R_+^2}\times \Xi_{\R}\times \Xi_{\R} \times \Xi_{\mathcal Z}\times\Xi_{\bbU}\times\R\times\calR^\bbG\times\calR^\bbG\times\mathcal{I}^\bbG\times\Omega\longrightarrow\R.
 \end{equation*}

The Hamiltonians functionals are defined as
\begin{align}\label{Hamilton}
\calH_i(t,\lambda,x,y_i,\theta_i,u, y_i^0, \zeta_i,p_i,q_i) &:=H_0^i(t,\lambda,x,y_i,\theta_i,u,y_i^0,\zeta_i,p_i,q_i)\\
&\quad +H_1^i(t,\lambda,x,y_i,\theta_i,u,y_i^0,\zeta_i,p_i,q_i)\nonumber,
\end{align}
where
\begin{align*}
	&H_0^i(t,\lambda,x,y_i,\theta_i,u,y_i^0,\zeta_i,p_i,q_i)\\
	&\quad:=F_i(t,u(t),x(t),y_i(t))+b(t,t,\lambda_t,u(t),x(t))p_i(t)+\kappa(t,t,0,\lambda_t,u(t),x(t))q_i(t,0)\lambda_t^B\\
	&\quad+\int_{\R_0}\kappa(t,t,z,\lambda_t,u(t),x(t))q_i(t,z)\lambda_t^H\nu(dz)+g_i(t,\lambda_t,u(t),x(t),y_i(t),\theta_i)\zeta_i(t)\\
	&H_1^i(t,\lambda,x,y_i,\theta_i,u,y_i^0,\zeta_i,p_i,q_i)\\
	\quad &:=\int_0^t\partial_t b(t,s,\lambda_s,u(s),x(s))ds\ p_i(t)+\int_0^t\!\!\!\int_{\R}\partial_t\kappa(t,s,z,\lambda_s,u(s),x(s))\calD_{s,z}p_i(t)\Lambda(dsdz),
\end{align*}
where $(p_i, q_i)$ and $\zeta_i$ are the solution to a forward-backward system of adjoint equation which will be introduced in \eqref{dp_i}, \eqref{dZ_i} here below.
For notational simplicity, from now we will write
\begin{align*}
   \calH_i(t) := \calH_i(t,\lambda, X, Y_i, \Theta_i(t,\cdot), u, y_i^0, \zeta_i,p_i,q_i),
\end{align*}
whenever we refer to the evaluation of $\calH_i(t,\lambda,x,y_i,\theta_i,u, y_i^0, \zeta_i,p_i,q_i)$ in $\t$, $\lambda \in \calL$, $X$ in \eqref{X_t}, $(Y_i,\Theta_i, M_i)$ in \eqref{Y_t}, $u\in\calA^{\caltildeE}_1\times \calA^{\caltildeE}_2$, $(p_i,q_i)$ in \eqref{dp_i} and $\zeta_i$ in \eqref{dZ_i}, and analogously for $H^i_0(t)$, $H^i_1(t)$.

For $\theta \in \mathcal{Z}$, we denote with $\partial_{\theta_{0}}$ the partial derivative with respect to $\theta(0)$ and $\nabla_{\theta_{z}}$ denotes the Fréchet derivative with respect to $\theta(z)$, $z\neq 0$. We also denote with $\frac{d}{d\nu}\nabla_{\theta_{z}}\calH_i(t)$ the Radon-Nikodim derivative of $\nabla_{\theta_{z}}\calH_i(t)$ with respect to $\nu(dz)$. \\
The couple $(p_i,q_i)$ in \eqref{Hamilton} satisfies the stochastic backward equation
\begin{equation}\label{dp_i}
	\begin{cases}
		dp_i(t)&=-\partial_x\calH_i(t)dt+\int_\R q_i(t,z)\mu(dtdz),\quad \t,\\
		p_i(T)&=\partial_x\varphi_i(X(T))+h_i(X(T))\zeta_i(T),
	\end{cases}	
\end{equation}
and $\zeta_i$ satisfies stochastic forward equation
\begin{equation}\label{dZ_i}
	\begin{cases}
		d\zeta_i(t)&=\partial_y\calH_i(t)dt+\partial_{\theta_{0}}\calH_i(t)dB(t)+\int_{\R_0}\frac{d}{d\nu}\nabla_{\theta_{z}}\calH_i(t)\tH(dtdz),\quad \t, \\
		\zeta_i(0)&=\partial_y\psi_i(y_i^0), 
	\end{cases}
\end{equation}
Also above, $\calD_{s,z}p_i (t)$ denotes the NA-derivative of $p_i(t)$ as in Theorem \ref{C-Ok-nonanticipating}. 

From now on we make the following crucial assumptions for $i=1,2$:
\begin{itemize}
    \item For all $i=1,2$, $h_i(X(T))\zeta_i(T)\in L^2(dP)$ (sufficient condition for this is to have $h_i$ bounded).
	\item $\calH_i(t)$ is well defined, Fréchet differentiable with respect to $x, y_i, \theta_i(0)$ $\theta_i(z)$, $z\neq 0$.
	\item $\frac{d}{d\nu}\nabla_{\theta_{z}}\calH_i(t)$ is well defined.
\end{itemize}

\begin{obs}
 Notice that the system \eqref{dp_i}-\eqref{dZ_i} is a partially coupled system of non Volterra forward-backward stochastic differential equations. Sufficient conditions that guarantee the existence of a solution of $\zeta_i \in L^2(dt\times dP)$, $(p_i,q_i)\in L^2(dt\times dP)\times\mathcal{I}^{\bbG}$ are, for $i=1,2$:
	\begin{enumerate}
	    \item $\E[\int_0^T\partial_x\calH_i(t,\lambda,x,y,\theta, u, 0,0,0)^2dt]<\infty$,
	    \item $\partial_x\calH_i(t)$ is Lipschitz with respect to $p$ and $q$ uniformly with respect to the other variables,
	    \item $\int_0^T\partial_{\theta_0}\calH_i(t)\lambda_t^Bdt<\infty$ $\mathbb{P}-a.s.$  and $\int_0^T\!\!\!\int_{\R_0}\frac{d}{d\nu}\nabla_{\theta_z}\calH_i(t,z)\lambda_t^H\nu(dz)dt<\infty$ $\mathbb{P}$-a.s.,
	    \item $\partial_y\calH_i(t), \partial_{\theta_0}\calH_i(t)$ and $\frac{d}{d\nu}\nabla_{\theta_z}\calH_i(t,z)$ are Lipschitz with respect to $x$, uniformly with respect to the other variables.
	\end{enumerate}

	Being the system only partially coupled, in fact, we can find $\zeta_i$ that solves the forward equation \eqref{dZ_i}, substitute it in \eqref{dp_i} and then find the pair $(p_i,q_i)$ solution to the backward equation \eqref{dp_i}. For more details on the solution of both the forward and the backward SDE of this type we refer to \cite{DiNunno1}. 
\end{obs}
\begin{obs}
    Following up Remark \ref{remark funzionale}, we notice that, when considering the forward equation \eqref{X_t} as a functional SDE \eqref{dX_t}, we would still obtain the Hamiltonian functional \eqref{Hamilton} (this is in line with e.g. \cite{DhalGames}).
\end{obs}

\begin{notz}\label{NOTZ1}
		From now on whenever it is clear the dependence of the coefficients from the processes $X$, $Y$, $\lambda$, $u$,..., we will use the simplified notation $b(t,s)=b(t,s,\lambda_s,u(s),X(s))$, $\hat b(t,s)=b(t,s,\lambda_s,\hatu(s),\hatX(s))$, $\kappa(t,s,z)=\kappa(t,s,z,\lambda_s,u(s),X(s))$, and also $\hat\kappa(t,s,z)=$
		$\hat\kappa(t,s,z,\lambda_s,\hatu(s),\hatX(s))$. We proceed similarly for $F_i$, $\hat F_i$ and $g_i$, $\hat g_i$.
\end{notz}

\subsection{A sufficient maximum principle}

The Nash equilibrium of the stochastic game is searched within the $\bbF$-adapted controls. In line with \cite{DiNunno1,DiNunnoGiordano}, for what concerns the interplay of the various information flows, and of \cite{Oksendal1}, for what concerns the treatment of Volterra structures, we introduce the $\bbF$-adapted processes $\calH^{\bbF}_i(t,\lambda, x, y_i, \theta_i,u,y_i^0, \zeta_i, p_i, q_i )$, $t\in [0,T]$, $\lambda\in \Xi_{\R^2_+}$, $x, y_i\in \Xi_\R$, $\theta_i\in \Xi_{\mathcal{Z}}$, $u\in \Xi_\R$, $y_i^0\in \R$, $\zeta_i, p_i \in \mathcal{R}^\bbG$, $q_i\in \mathcal{I}^\bbG$, as follows:
	\begin{align}\label{HamiltonF}
&\calH^{\bbF}_i(t,\lambda,x,y_i,\theta_i,u,y_i^0, \zeta_i,p_i,q_i)\\
&:=\E\left[\calH_i(t,\lambda,x,y_i,\theta_i,u, y_i^0, \zeta_i,p_i,q_i)|\F_t\right]\nonumber\\
&=H^{\bbF,i}_0(t,\lambda,x,y_i,\theta_i,u,y_i^0, \zeta_i,p_i,q_i)+H^{\bbF,i}_1(t,\lambda,x,y_i,\theta_i,u,y_i^0, \zeta_i,p_i,q_i) \nonumber,
	\end{align}
where
\begin{align*}
&H^{\bbF,i}_0(t,\lambda,x,y_i,\theta_i,u,y_i^0, \zeta_i,p_i,q_i)\\
&\quad:=F_i(t,u,x(t),y_i(t))+b(t,t,\lambda_t,u(t),x(t))\E[p_i(t)|\F_t]\\
&\qquad+\kappa(t,t,0,\lambda_t,u(t),x(t))\E[q_i(t,0)|\F_t]\lambda_t^B\\
&\qquad+\int_{\R_0}\kappa(t,t,z,\lambda_t,u(t),x(t))\E[q_i(t,z)|\F_t]\lambda_t^H\nu(dz)+g_i(t,\lambda_t,u(t),x(t),y_i(t),\theta_i)\E[\zeta_i|\F_t]\\
&H^{\bbF,i}_1(t,\lambda,x,y_i,\theta_i,u,y_i^0, \zeta_i,p_i,q_i)\\
&\quad:=
\int_0^t\partial_t b(t,s,\lambda_s,u(s),x(s))ds\ \E[p_i(t)|\F_t]\\&\qquad+\int_0^t\!\!\!\int_{\R}\partial_t\kappa(t,s,z,\lambda_s,u(s),x(s))\E[\calD_{s,z}p_i(t)|\F_t]\Lambda(dsdz)
\end{align*}
With notation analogous to the case with information flow $\bbG$, we will write $\calH_i^{\bbF}(t)$ to refer to the functional \eqref{HamiltonF} evaluated in $\t$, $\lambda \in \calL$, $X$ in \eqref{X_t}, $(Y_i,\Theta_i(t,\cdot), M_i)$ in \eqref{Y_t}, $u\in\calA^{\calE}_1\times \calA^{\calE}_2$, $(p_i,q_i)$ in \eqref{dp_i} and $\zeta_i$ in \eqref{dZ_i}. Analogous notation is used for $H_0^{\bbF,i}(t)$ and $H_1^{\bbF,i}(t)$. We have the following:

\begin{teo}\label{SufficientMaximumTheoremF}
	Let $\hat u= (\hatu_1,\hatu_2)\in\calA^{\calE}_1\times\calA^{\calE}_2$ and assume that the corresponding solutions $\hatX$, $(\hatY_i, \hat\Theta_i, \hat M_i)$, $(\hatp_i, \hatq_i)$, $\hat\zeta_i$ of equations \eqref{X_t}, \eqref{Y_t}, \eqref{dp_i}, and \eqref{dZ_i} (with $y_i^0 = \hatY_i(0)$ in the initial value) exist for $i=1,2$. 
	We assume that the functionals $h_i$ in \eqref{Y_t} are concave. We also consider performance functionals \eqref{J} with concave $\varphi_i$ and $\psi_i$.
	Assume that $\calH^\bbF_i$ \eqref{HamiltonF} satisfy the following conditions:
	\begin{enumerate}
	\item[(i)]
	The maps
		$$x,y,\theta \longmapsto \calH_i^\bbF(t,\lambda,x,y,\theta,u,y^0,\zeta,p,q)$$
		are concave for all $\t$, $\lambda \in \calL$, $u\in\bbU$, $y^0\in \R$, $z,p \in\calR^\bbG$, $q\in\mathcal I^\bbG$.\\
	\item[(ii)] The following extremes are achieved:
		\begin{align*}
		\sup_{v\in\calA^{\calE}_1}&\E\left[\calH^{\bbF}_1(t,\lambda,\hatX,\hatY_1,\hat\Theta_1(t,\cdot),v,\hatu_2,\hat \zeta_1,\hatp_1,\hatq_1)|\calE_t^{(1)}\right]\\
		&=\E\left[\calH^{\bbF}_1(t,\lambda,\hatX,\hatY_1,\hat\Theta_1(t,\cdot),\hatu_1,\hatu_2,\hat \zeta_1,\hatp_1,\hatq_1)|\calE_t^{(1)}\right],
		\end{align*}
		and
		\begin{align*}
		\sup_{v\in\calA^{\calE}_2}&\E\left[\calH^{\bbF}_2(t,\lambda,\hatX,\hatY_2,\hat\Theta_2(t,\cdot),\hatu_1,v,\hat \zeta_2,\hatp_2,\hatq_2)|\calE_t^{(2)}\right]\\
		&=\E\left[\calH^{\bbF}_2(t,\lambda,\hatX,\hatY_2,\hat\Theta_2(t,\cdot),\hatu_1,\hatu_2,\hat \zeta_2,\hatp_2,\hatq_2)|\calE_t^{(2)}\right].
		\end{align*}
		\end{enumerate}
		Then $\hatu=(\hatu_1,\hatu_2)$ is a Nash equilibrium with respect to the information flow $\bbF$ for the stochastic game \eqref{X_t},\eqref{Y_t},\eqref{J}.
\end{teo}

\begin{proof}
	We consider $J_1$, we prove that $J_1(u_1,\hatu_2)\leq J_1(\hatu_1,\hatu_2)$, for all $u_1\in \calA^{\calE}_1$ and that $J_1(\hatu_1,u_2)\leq J_1(\hatu_1,\hatu_2)$, for all $u_2\in \calA^{\calE}_2$. We can then follow the same steps for $J_2$ and conclude. 
	We start by the first inequality. 
	
	Let $u_1\in\calA^{\calE}_1$ and denote $X, (Y_1,\Theta_1,M_1)$ the solutions of \eqref{X_t}, \eqref{Y_t} with $u=(u_1, \hat u_2)$ and $\hat X, (\hat Y_1,\hat \Theta_1,\hat M_1)$ those with $\hat u=(\hat u_1, \hat u_2)$.  Also let $(\hatp_1,\hatq_1,\hat\zeta_1)$ to be the solutions to the adjoint systems  \eqref{dp_i}-\eqref{dZ_i} with $\hatu$ and initial value with $y_i^0=\hatY_1(0)$. Correspondingly, we shall consider \eqref{HamiltonF} with
	\begin{align*}
		\calH_1^{\bbF}(t)&:=\calH_1^{\bbF} (t,\lambda,X,Y_1,\Theta_1(t,\cdot),(u_1,\hatu_2), Y_1(0), \hat \zeta_1,\hatp_1,\hatq_1),\\
		\hat\calH_1^{\bbF}(t)&:=\hat\calH_1^{\bbF}(t,\lambda,\hatX,\hatY_1,\hat\Theta_1(t,\cdot),(\hatu_1,\hatu_2), \hatY_1(0), \hat \zeta_1,\hatp_1,\hatq_1).
	\end{align*}
 Correspondingly, we have $H_0^{\bbF, i}(t)$, $\hat H_0^{\bbF, i}(t)$, $H_1^{\bbF, i}(t)$, $\hat H_1^{\bbF, i}(t)$.
We study the difference
	\begin{equation*}
	\Delta_1:=J_1(u_1,\hatu_2)-J_1(\hatu_1,\hatu_2)=I_1+I_2+I_3,
	\end{equation*}
	as sum of three elements:
	\begin{align}
	I_1&:=\E\left[\int_0^T\left(F_1(t,u(t),X(t),Y_1(t))-F_1(t,\hatu(t),\hatX(t),\hatY_1(t) \right)dt\right],\label{I1}\\
	I_2&:=\E[\varphi_1(X(T))-\varphi_1(\hatX(T))],\label{I2}\\
	I_3&:=\E[\psi_1(Y_1(0))-\psi_1(\hatY_1(0))]\label{I3}.
	\end{align}
Note that by introducing a sequence of stopping times similarly to what has been done in \cite{OksendalGames} Theorem 6.4, we can assume that all the local martingales appearing in the upcoming computations are martingales. In particular, the expectations of the terms with $\mu(dtdz)$-integrals and $dM_i$-integrals are all 0. 

Starting from $I_1$ we see that
	\begin{align*}
	I_1&=\E\Bigg[\int_0^T\Big\{ H_0^{\bbF,1}(t)-\hat{H}_0^{\bbF,1}(t)-\E\left[\hatp_1(t)|\F_t\right]\left(b(t,t)-\hatb(t,t)\right)\Big\}dt \\
	&\quad-\int_0^T\int_{\R}\E\left[\hatq_1(t,z)|\F_t\right]\left(\kappa(t,t,0)-\hatk(t,t,z)\right)\Lambda(dtdz)\\
	&\quad-\int_0^T\E\left[\hat\zeta_1(t)|\F_t\right]\left(g_1(t)-\hatg_1(t)\right) dt\Bigg].
	\end{align*}
To study $I_2$, we use the concavity condition of $\varphi_1$, the Transformation Rule (\cite{DiNunnoGiordano} Lemma 3.4) applied to $X(t)$ and Itô's formula for the product. Hence we obtain
	\begin{align*}
	I_2&\leq\E\left[\partial_x\varphi(\hatX(T))(X(T)-\hatX(T))\right]\\
	&=\E\left[\hatp_1(T)(X(T)-\hatX(T))\right]
	-\E\left[\hat\zeta_1(T) \partial_x h_1(\hatX(T))(X(T)-\hatX(T))\right]\\
	&= \E\left[\int_0^T\left\{\hatp_1(t)\left(b(t,t)-\hatb(t,t)\right)+\int_0^t\left(\partial_t b(t,s)-\partial_t \hatb(t,s)\right)ds\ \hatp_1(t)\right\} dt \right.\\
	&\quad+\int_0^T\int_0^t\!\!\!\int_{\R}\Big(\partial_t \kappa(t,s,z)-\partial_t \hatk(t,s,z)\Big)\mu(dsdz)\hat p_1 (t)dt \\
	&\quad\left. - \int_0^T \partial_x\hat{\calH}^{\bbF}_1(t)\Big(X(t)-\hatX(t)\Big) dt+\int_0^T\!\!\!\int_{\R} \hat q_1(t,z)\Big(\kappa(t,t,z)-\hatk(t,t,z)\Big)\Lambda(dtdz)\right]\\
	&\quad-\E\left[\hat\zeta_1(T) \partial_x h_1(\hatX(T))(X(T)-\hatX(T))\right].
	\end{align*}
	Exploiting the Fubini's theorem together with the duality formula of Corollary \ref{Yab} we have that
    \begin{align*}
	\E\left[\int_0^T\left(\int_0^t\!\!\!\int_{\R}\partial_t \kappa(t,s,z)\mu(dsdz)\right)\hatp_1(t) dt\right]
	&=\int_0^T\E\left[\left(\int_0^t\!\!\!\int_{\R}\partial_t \kappa(t,s,z)\mu(dsdz)\right)\hatp_1(t) \right]dt\\
	&=\int_0^T\E\left[\!\!\!\int_0^t\!\!\!\int_\R\partial_t\kappa(t,s,z)\calD_{s,z}\hat p_1(t)\Lambda(dsdz) \right]dt\\
	&=\E\left[\int_0^T\!\!\!\int_0^t\!\!\!\int_\R\partial_t\kappa(t,s,z)\calD_{s,z}\hat p_1(t)\Lambda(dsdz)dt \right]
	\end{align*}
	Using now the concavity of $\psi_1$ and $h_1$, together with the Itô formula, we get that
	\begin{align}
	I_3&\leq\E\Big[\partial_y\psi_1(\hatY(0))\left(Y_1(0)-\hatY_1(0)\right)\Big]\label{I3_inequality}\\
	&=\E\left[\hat\zeta_1(0)(Y_1(0)-\hatY_1(0))\right]\nonumber\\
	&=\E\left[\hat\zeta_1(T)(Y_1(T)-\hatY_1(T))\right]\nonumber\\
	&
	-\E\left[\int_0^T\left(Y_1(t)-\hatY_1(t)\right)d\hat\zeta_1(t)+\int_0^T\zeta_1(t)\left(dY_1(t)-d\hatY_1(t)\right)\right.\nonumber\\
	& \left.
	+\hspace{-1mm}\int_0^T \hspace{-2mm}\Big\{ \partial_{\theta_{0}} \hat{\calH}^{\bbF}_1(t)\left(\Theta_1(t,0)-\hat\Theta_1(t,0)\right)\lambda^B_t 
	\hspace{-1mm}+ \hspace{-1mm}\int_{\R_0}\frac{d}{d\nu}\nabla_{\theta_{z}}\hat{\calH}^{\bbF}_1(t)\left(\Theta_1(t,z)-\hat\Theta_1(t,z) \right)\lambda^H_t \nu(dz) \Big\}dt \right]\nonumber\\
	&\leq\E\left[\hat\zeta_1(T)\partial_x h_1(\hatX(T))(X(T)-\hatX(T))\right]
	-\E\left[\int_0^T \Big\{ \partial_y\hat{\calH}^{\bbF}_1(t)\left(Y_1(t)-\hatY_1(t)\right)\right.\nonumber\\
	&+ \partial_{\theta_{0}}\hat{\calH}^{\bbF}_1(t)\left(\Theta_1(t,0)-\hat\Theta_1(t,0)\right)\lambda^B_t 
	+\int_{\R_0}\frac{d}{d\nu}\nabla_{\theta_{z}}\hat{\calH}^{\bbF}_1(t)\left(\Theta_1(t,z)-\hat\Theta_1(t,z)\right) \nu(dz)\lambda^H_t\nonumber\\
	&\left.+\hat\zeta_1(t)\left(g(t)-\hatg(t)\right) \Big\} dt\right] \nonumber.
	\end{align}
	Hence, summing up the three estimates above, conditioning on $\F_t$ and recalling \eqref{HamiltonF}, we get that
	\begin{align}
	\Delta&=I_1+I_2+I_3\nonumber \\
	&\leq \E\Bigg[\int_0^T\Big\{\calH^{\bbF}_1(t)-\hat{\calH}^{\bbF}_1(t)-\partial_x \hat{\calH}^{\bbF}_1(t) \left(X(t)-\hatX(t)\right)\nonumber\\
	&\quad-\partial_y\hat{\calH}^{\bbF}_1(t)\left(Y_1(t)-\hatY_1(t)\right)-\partial_{\theta_{0}}\hat{\calH}^{\bbF}_1(t)\left(\Theta_1(t,0)-\hat\Theta_1(t,0)\right) \nonumber \\
	&\quad-\int_{\R_0}\frac{d}{d\nu}\nabla_{\theta_{z}}\hat{\calH}^{\bbF}_1(t) \left(\Theta_1(t,z)-\hat\Theta_1(t,z)\right) \nu(dz)\Big\}dt\Bigg],\nonumber
	\end{align}
is non-negative by concavity of $\calH^{\bbF}_1$, assumption (i). Thanks to (ii), after taking the conditional expectation under $\calE^{(1)}_t$, we obtain that, for all $u_1\in\calA^{\calE}_1$,
	\begin{equation*}\label{Nash1}
	J_1(u_1,\hatu_2)\leq J_1(\hatu_1,\hatu_2).
	\end{equation*}
	The inverse inequality can be proved by analogous arguments.
\end{proof}
\noindent Similarly to what we did for the information $\E^{(i)}$, we can also formulate a sufficient maximum principle with respect to the enlarged information $\tilde\E^{(i)}\subset \bbG$, see Remark \ref{EnlargedE}. Even though the result is possibly not directly applicable, given the nature of the enlarged filtration $\bbG$, the study has mathematical validity. Recall the compact notation introduced in Notation \ref{NOTZ1}.

\begin{prop}\label{SufficientMaximumTheorem}
	Let $\hatu= (\hatu_1,\hatu_2)\in\calA^{\caltildeE}_1\times\calA^{\caltildeE}_2$ and assume that the corresponding solutions $\hatX$, $(\hatY_i,\hat\Theta_i)$,$(\hatp_i, \hatq_i)$, $\hat\zeta_i$ of equations \eqref{X_t}, \eqref{Y_t},\eqref{dp_i} and \eqref{dZ_i} (with $y_i^0=\hatY_i(0)$ in the initial value) exist for $i=1,2$. Here we assume that the functionals $h_i$ in \eqref{Y_t} are concave. We also consider performance functionals \eqref{J} with concave $\varphi_i$ and $\psi_i$.
	Assume that $\calH_i$ \eqref{Hamilton} satisfy the following conditions:
	\begin{enumerate}
	\item[(i)]
		The maps
		$$x,y,\theta\longmapsto \calH_i(t,\lambda,x,y,\theta,y^0,\zeta,p,q)$$
		are concave for all $\t$, $\lambda \in \calL$, $y^0\in\R$, $z,p \in\calR^\bbG$, $q\in\mathcal I^\bbG$.\\
\item[(ii)] The extremes are achieved
		\begin{align*}
		\sup_{v\in\calA^{\caltildeE}_1}&\E\left[\calH_1(t,\lambda,\hatX,\hatY_1,\hat\Theta_1(t,\cdot),v,\hatu_2,\hatY_1(0),\hat \zeta_1,\hatp_1,\hatq_1|\caltildeE_t^{(1)}\right]\\
		&=\E\left[\calH_1(t,\lambda,\hatX,\hatY_1,\hat\Theta_1(t,\cdot),\hatu_1,\hatu_2,\hatY_1(0), \hat \zeta_1,\hatp_1,\hatq_1)|\caltildeE_t^{(1)}\right],
		\end{align*}
		and
		\begin{align*}
		\sup_{v\in\calA^{\caltildeE}_2}&\E\left[\calH_2(t,\lambda,\hatX,\hatY_2,\hat\Theta_2(t,\cdot),\hatu_1,v,\hatY_2(0), \hat \zeta_2,\hatp_2,\hatq_2)|\caltildeE_t^{(2)}\right]\\
		&=\E\left[\calH_2(t,\lambda,\hatX,\hatY_2,\hat\Theta_2(t,\cdot),\hatu_1,\hatu_2, \hatY_2(0), \hat \zeta_2,\hatp_2,\hatq_2)|\caltildeE_t^{(2)}\right],
		\end{align*}
		\end{enumerate}
Then $\hatu=(\hatu_1,\hatu_2)$ is a Nash equilibrium for \eqref{X_t}-\eqref{Y_t}-\eqref{J} with respect to the information $\bbG$.
\end{prop}

\begin{proof}
	The arguments of Theorem \ref{SufficientMaximumTheorem} leading to 
	\begin{align*}
	\Delta&=I_1+I_2+I_3\\
	&\leq \E\left[\int_0^T\left\{\calH_1(t)-\hat{\calH}_1(t)-\partial_x \hat{\calH}_1(t) \left(X(t)-\hatX(t)\right)\right.\right.\\
	&\quad-\partial_y\hat{\calH}_1(t)\left(Y_1(t)-\hatY_1(t)\right)-\partial_{\theta_{0}}\hat{\calH}_1(t)\left(\Theta_1(t,0)-\hat\Theta_1(t,0)\right)\\
	&\quad\left.\left.-\int_{\R_0}\frac{d}{d\nu}\nabla_{\theta_{z}}\hat{\calH}_1(t) \left(\Theta_1(t,z)-\hat\Theta_1(t,z)\right)\nu(dz)\right\}dt\right],
	\end{align*}
hold. Then the conditioning with respect to $\tilde\E^{(i)}$ as a sub-filtration of $\bbG$ is taken. 
\end{proof}

\begin{obs}
    Notice that Theorem \ref{SufficientMaximumTheoremF} and Proposition \ref{SufficientMaximumTheorem} still hold even without the hypothesis necessary for the Transformation rule to hold (see \eqref{HPkap}) whenever it is possible to apply an Itô formula for the product  $\mathbb E[\hat p(T)(X(T) - \hat X(T)]$. This is the case, e.g. when $\kappa$ is not of Volterra-type.
\end{obs}

\subsection{The zero-sum game case} 

In the zero-sum game, we recall that a player's gain corresponds to the other player's loss. This translates directly into having
\begin{equation}\label{J(u) zero sum}
	J_1(u_1,u_2)=-J_2(u_1,u_2).
\end{equation}
In this case one can regard the Nash equilibrium of the game $\hatu=(\hatu_1,\hatu_2)$ as a saddle point for the functional
\begin{equation}\label{J(u) zero sum 2}
    J(u_1,u_2):=J_1(u_1,u_2),
\end{equation}
$(u_1,u_2)\in \calA_1^{\cdot}\times\calA_2^\cdot$, where the controls may refer to the information flows $\E^i\subseteq \bbF$ or $\tilde{\E}^i\subseteq\bbG$.

To prove this statement it is sufficient to notice that from Definition \ref{Nash} and Remark \ref{EnlargedE} we have that 
\begin{equation*}
	J_1(u_1,\hatu_2)\leq J_1(\hatu_1,\hatu_2)=-J_2(\hatu_1,\hatu_2)\leq -J_2(\hatu_1,u_2),
\end{equation*}
and hence, for all $u_1,u_2$,
\begin{equation*}
	J(u_1,\hatu_2)\leq J(\hatu_1,\hatu_2)\leq J(\hatu_1,u_2).
\end{equation*}
As a consequence one has that
\begin{align*}
	\inf_{u_2\in\calA^{\cdot}_2}\sup_{u_1\in\calA^{\cdot}_1} J(u_1,u_2)\leq J(\hatu_1,\hatu_2)\leq \sup_{u_1\in\calA^{\cdot}_1}\inf_{u_2\in\calA^{\cdot}_2} J(u_1,u_2),
\end{align*}
 and since it is always true that $\inf\sup f\geq \sup \inf f $ for a given function $f$, we get the desired result:
\begin{align}\label{saddle}
\inf_{u_2\in\calA^{\cdot}_2}\sup_{u_1\in\calA^{\cdot}_1} J(u_1,u_2)= J(\hatu_1,\hatu_2)= \sup_{u_1\in\calA^{\cdot}_1}\inf_{u_2\in\calA^{\cdot}_2} J(u_1,u_2).
\end{align}
The study of a zero-sum game \eqref{J(u) zero sum}-\eqref{saddle} involves only one Hamiltonian and only one set of adjoint equations. In fact, by defining
\begin{equation}\label{definition coefficients zero}
    g_1=g_2=:g, \ h_1=h_2=:h, \ F_1=-F_2=:F, \ \varphi_1 = -\varphi_2 =:\varphi, \ \text{and} \ \psi_1=\psi_2=:\psi,
\end{equation}
we see that the Hamiltonian functionals \eqref{Hamilton} become
\begin{align*}
	\calH_1(t)&=F(t,u(t),X(t),Y(t))+b(t,t,\lambda_t,u(t),X(t))p_1+\kappa(t,t,0,\lambda_t,u(t),X(t))q_1(t,0)\lambda_t^B\\
	&\quad+\int_{\R_0}\kappa(t,t,z,\lambda_t,u(t),X(t))q_1(t,z)\lambda_t^H\nu(dz)+g(t,\lambda_t,u(t),X(t),Y(t),\Theta(t,\cdot))\zeta_1(t)\\
	&\quad+\int_0^t\partial_t b(t,s,\lambda_s,u(s),X(s))ds\ p_1(t) \\
	&\quad+\int_0^t\!\!\!\int_{\R}\partial_t\kappa(t,s,z,\lambda_s,u(s),X(s))\calD_{s,z}p_1(t)\Lambda(dsdz),
\end{align*}
and 
\begin{align*}
\calH_2(t)&=F(t,u(t),X(t),Y(t))+b(t,t,\lambda_t,u(t),X(t))p_2+\kappa(t,t,0,\lambda_t,u(t),X(t))q_2(t,0)\lambda_t^B\\
&\quad+\int_{\R_0}\kappa(t,t,z,\lambda_t,u(t),X(t))q_2(t,z)\lambda_t^H\nu(dz)+g(t,\lambda_t,u(t),X(t),Y(t),\Theta(t,\cdot))\zeta_2(t)\\
&\quad+\int_0^t\partial_t b(t,s,\lambda_s,u(s),X(s))ds\ p_2(t)\\ &\quad+\int_0^t\!\!\!\int_{\R}\partial_t\kappa(t,s,z,\lambda_s,u(s),X(s))\calD_{s,z}p_2(t)\Lambda(dsdz),
\end{align*}
 where $u=(u_1,u_2)\in\calA_1\times\calA_2$, $X$ is defined in \eqref{X_t} and $(Y,\Theta, M)$ is given in $\eqref{Y_t}$ with \eqref{definition coefficients zero} applied. From \eqref{dp_i} and \eqref{dZ_i}, with \eqref{definition coefficients zero}, we have
\begin{equation}\label{adjoint_zero_sum}
    \zeta_2(t)=-\zeta_1(t),\quad \text{ and } \quad (p_2(t),q_2(t,\cdot))=-(p_1(t),q_1(t,\cdot)),\quad \t.
\end{equation}
  Then we can write
\begin{align}\label{adjoint_zero}
	-\calH_2&(t,\lambda,X,Y,\Theta(t,\cdot),u_1,u_2,\zeta_2,p_2,q_2)\nonumber\\
	&=\calH_1(t,\lambda,X,Y,\Theta(t,\cdot),u_1,u_2,\zeta_1,p_1,q_1) .
\end{align}
Hence in the zero-sum game case, we also need only one Hamiltonian functional of type \eqref{Hamilton} and one triple of adjoint processes. We denote such Hamiltonian functional
\begin{equation}\label{HamiltonZero}
	\calH(t):=-\calH_2(t)=\calH_1(t),
\end{equation}
and we define
\begin{align*}
	(\zeta,p,q):=(\zeta_1,p_1,q_1)=(-\zeta_2,-p_2,-q_2), \ \text{ for all }\t.
\end{align*}
Also, we rewrite the performance functional \eqref{J(u) zero sum 2} as
\begin{equation}\label{PerformanceZeroSum}
	J(u_1,u_2)=\E\left[\int_0^T F(t,u(t),X(t),Y(t))dt+\varphi(X(T))+\psi(Y(0))\right].
\end{equation}
Analogously, when dealing with the information flow $\bbF$, we can consider only one Hamiltonian functional of type \eqref{HamiltonF}
\begin{equation}\label{HamiltonZeroF}
\calH^{\bbF}:=-\calH^{\bbF}_2=\calH^{\bbF}_1.
\end{equation}
We can thus reformulate Theorems \ref{SufficientMaximumTheoremF} for the zero-sum game cases:

\begin{teo}\label{SufficientMaximumZeroSumTheoremF}
Let $\hatu= (\hatu_1,\hatu_2)\in\calA^{\calE}_1\times\calA^{\calE}_2$ and assume that the corresponding solution $\hatX$ of \eqref{X_t} and $(\hatY,\hat\Theta)$,$(\hatp, \hatq)$, $\hat\zeta$ corresponding to \eqref{Y_t},\eqref{dp_i} and \eqref{dZ_i} with \eqref{definition coefficients zero} in the zero-sum context exist. Here the functional $h$ is assumed to be concave. We also consider the performance functional \eqref{PerformanceZeroSum} with concave $\varphi$ and $\psi$.
	Assume that $\calH^{\bbF}$ \eqref{HamiltonZeroF} satisfies the following conditions:
	\begin{enumerate}
	\item[(i)]
		The maps
		$$x,y,\theta\longmapsto \calH^{\bbF}(t,\lambda,x,y,\theta,y^0,\zeta,p,q)$$
		are concave for all $\t$, $\lambda \in \calL$, $y^0\in\R$, $z,p \in\calR^\bbG$, $q\in\mathcal I^\bbG$.\\
\item[(ii)] The extremes are achieved
		\begin{align*}
		\sup_{v\in\calA^{\calE}_1}&\E\left[\calH^{\bbF}(t,\lambda,\hatX,\hatY,\hat\Theta(t,\cdot),v,\hatu_2,\hatY(0),\hat \zeta,\hatp,\hatq)|\calE_t^{(1)}\right]\\
		&=\E\left[\calH^{\bbF}(t,\lambda,\hatX,\hatY,\hat\Theta(t,\cdot),\hatu,\hatu_2,\hatY(0), \hat \zeta,\hatp,\hatq)|\calE_t^{(1)}\right],
		\end{align*}
		and similarly
		\begin{align*}
		\sup_{v\in\calA^{\caltildeE}_2}&\E\left[\calH^{\bbF}(t,\lambda,\hatX,\hatY,\hat\Theta(t,\cdot),\hatu_1,v,\hatY(0), \hat \zeta,\hatp,\hatq)|\calE_t^{(2)}\right]\\
		&=\E\left[\calH^{\bbF}(t,\lambda,\hatX,\hatY,\hat\Theta(t,\cdot),\hatu_1,\hatu_2, \hatY(0), \hat \zeta,\hatp,\hatq)|\calE_t^{(2)}\right],
		\end{align*}
		\end{enumerate}
Then $\hatu=(\hatu_1,\hatu_2)$ is a saddle point for $J$.
\end{teo}

\noindent An analogous result holds also for zero-sum games when the information flow is $\bbG$ in line with Proposition \ref{SufficientMaximumTheorem} under \eqref{J(u) zero sum}-\eqref{J(u) zero sum 2}, \eqref{definition coefficients zero}.


\section{Optimal consumption and recursive utility}\label{section4}
\subsection{About optimal recursive utility, an example}

Hereafter we present a game with optimal recursive utility consumption. This illustration is related to \cite{Oksendal1} and the references therein, with the main  difference that we do consider a time-changed Lévy process as noise in the dynamics and two players competing against each other.
In this case our forward dynamics \eqref{X_t} are given by a cash flow $X$ exposed to a consumption rate $c_1(t)+c_2(t)$ with $c(t)=(c_1(t),c_2(t))$, $t\geq 0$ and $c_i(t)\geq0$ for $i=1,2$:
 \begin{equation}\label{ForwardUtility}
 	X(t)=X_0 + \int_0^tX(s)\alpha(t,s)-c_1(s)-c_2(s)ds+\int_0^t\!\!\!\int_\R\pi(t,s,z)\mu(dsdz), \quad \t ,
 \end{equation}
 where we assume that $X_0\in\R$, $\alpha:[0,T]^2\longrightarrow \R$ and $\pi:[0,T]^2\times \R \longrightarrow \R$ are such that $\alpha,\pi$ bounded and with derivatives $\partial_t\alpha(t,s)$ and $\partial_t\pi(t,s,\cdot)$ in $L^2(dt)$. We also assume that $\partial_t\pi(t,s,z)$ is Lipschitz with respect to $\t$.
 Inspired by \cite{Duffie-Epstein}, for each player $i=1,2$ we intend the backward dynamics \eqref{Y_t} as recursive utility processes $Y_i$, $i=1,2$ given by
 \begin{align}\label{backwardUtility}
 	dY_i(t)=-\left[\gamma_i(t)Y_i(t)+\ln(c_i(t)X(t))\right]dt+\int_\R\Theta_i(t,z)\mu(dtdz)+dM_i(t), \quad \t,
 \end{align}
 with $Y(T)=0$. Above we assume that $\gamma_i$ is a $L^2(dt\times dP)$-integrable function.
The goal of each player $i=1,2$ is to maximize
\begin{equation}\label{optimalRecursiveF}
	J_i(c_i)=\E[Y_i(0)]
\end{equation}
	over the set of non-negative $\E^{(i)}$-predictable processes $c_i$, where $\E^{(i)}$ is the information available to each player.

We are interested in solving explicitly \eqref{ForwardUtility}. We start by defining inductively:
\begin{equation*}
    \begin{cases}
    X^{(0)}(t) &= X_0 - \int_0^t c_1(s)+c_2(s)ds + \int_0^t\int_\R\pi(t,s,z) \mu(dsdz) \\
    X^{(n+1)}(t) &= X_0 + \int_0^t \alpha(t,s) X^{(n)}(s) ds + \int_0^t c_1(s)+c_2(s)ds - \int_0^t\int_\R\pi(t,s,z) \mu(dsdz),
    \end{cases}
\end{equation*}
for $n\geq 0$, and
\begin{equation*}
    \begin{cases}
    A^{(1)}(t,s) &= \alpha(t,s)\\
    A^{(n+1)}(t,s) &= \int_s^t\alpha(t,v)A^{(n)}(v,s)dv, \quad n\geq 1.
    \end{cases}
\end{equation*}
Applying the Fubini theorem we can rewrite the $n$-th term above as
\begin{equation*}
    X^{(n)}(t) = \left(X_0-\int_0^tc_1(s)+c_2(s) ds + \int_0^t\!\!\!\int_\R\pi(t,s,z) \mu(dsdz)\right)\left(1+\int_0^t\sum_{k=0}^n A^{(k)}(t,s)ds\right).
\end{equation*}
Using a Picard iteration scheme we find that, by defining
\begin{equation}\label{M(t)}
\mathcal M(t) := 1+\int_0^t\sum_{k=1}^\infty A^{(k)}(t,s)ds, \quad \t,    
\end{equation}
the solution to \eqref{X_t Volterra} is given by
\begin{equation}\label{soluzione_esplicita_X}
    X(t) = \mathcal M(t) \left(X_0-\int_0^tc_1(s)+c_2(s) ds + \int_0^t\!\!\!\int_\R\pi(t,s,z) \mu(dsdz)\right),
\end{equation}
where we have that $\sum_{k=1}^\infty A^{(k)}(t,s)ds$ is converging for a.a. $(t,s)$ due to the boundedness of $\alpha$. In order to guarantee the well definition of $Y$ in \eqref{backwardUtility}, we will assume that $X(t)>0$ for all $\t$. To the best of our knowledge, currently there are no results that guarantee the positivity of a Volterra equation like \eqref{ForwardUtility}.
	
	Proceding as in Section 3.1, we consider now the Hamiltonian functionals \eqref{Hamilton}:
\begin{align}
	\calH_i(t)&=(\alpha(t,t)X(t)-c_1(t)-c_2(t))p_i(t)+\pi(t,t,0)q_i(t,0)\lambda_t^B\nonumber \\
	&\quad +\int_{\R_0}\pi(t,t,z)q_i(t,z)\lambda_t^H\nu(dz) +\int_0^t\partial_t\alpha(t,s) X(s) ds\  p_i(t)\nonumber\\
	&\quad+\int_0^t\!\!\!\int_\R\partial_t\pi(t,s,0) \calD_{s,z}p_i(t)\Lambda(ds,dz) +[\gamma_i(t)Y_i(t)+\ln(c_i)+\ln(X(t))]\zeta_i(t),\label{HamiltonExample}
\end{align}
and the corresponding backward system of adjoint equations \eqref{dp_i}-\eqref{dZ_i}:
\begin{align}
dp_i(t)&=-\left[\left(\alpha(t,t)+\int_0^t\partial_t\alpha(t,s)ds\right)p_i(t)+\frac{\zeta_i(t)}{X(t)}\right]dt+\int_\R q_i(t,z)\mu(dtdz),\nonumber\quad \t\\
p_i(T)&=0\label{dp(t)},
\end{align}
and
\begin{equation}
	d\zeta_i(t)=\gamma_i(t)\zeta_i(t)dt,\quad \label{dz(t)} \t,
\end{equation}
with initial condition $\zeta_i(0)=1$. The solution of \eqref{dz(t)} is immediate:
\begin{equation}\label{solution_z(t)}
	\zeta_i(t)=\exp\left\{-\int_0^t\gamma_i(s)ds\right\}.
\end{equation}
As for \eqref{dp(t)}, we know from \cite{DiNunno1} Theorem 5.1 that a sufficient condition that ensure the existence of a solution is that 
\begin{equation}\label{condition p}
    \E\left[\int_0^T\left(\frac{\zeta_i(t)}{X(t)}\right)^2 dt\right]<\infty,
\end{equation}
and in this case $p_i(t)$ has representation
\begin{equation}\label{soluzione_esplicita_p}
    p_i(t) =\E\left[\int_t^T\exp\left\{\int_t^s\left(\alpha(r,r)+\int_0^t\partial_r\alpha(r,v)dv\right) dr\right\}\frac{\zeta(s)}{X(s)}ds \right] 
\end{equation}
	If we now consider the problem under $(\bbF,\E^{(i)})$, our  functionals \eqref{HamiltonF} become
\begin{align}
		\calH^{\bbF}_i(t)&=(\alpha(t,t)X(t)-c_1(t)-c_2(t))\E[p_i(t)|\F_t]+\pi(t,t,0)\E[q_i(t,0)|\F_t]\lambda_t^B\nonumber\\
		&\quad+\int_{\R_0}\pi(t,t,z)\E[q_i(t,z)|\F_t]\lambda_t^H\nu(dz)+[\gamma_i(t)+\ln(c_i(t))+\ln(X(t))]\E[\zeta_i(t)|\F_t]\nonumber\\
		&\quad+\int_0^t\partial_t\alpha(t,s)X(s)ds\ E[p_i(t)|\F_t]+\int_0^t\!\!\!\int_{\R}\partial_t\pi(t,s,z)\E[\calD_{s,z}p_i(t)|\F_t]\Lambda(dsdz)\nonumber.
	\end{align}
Thanks to Theorem \ref{SufficientMaximumTheoremF}, a sufficient optimality condition is given by
\begin{equation}\label{Suff_cond}
	\partial_{c_i} \calH_i^{\bbF} =0.
\end{equation}
We thus claim that a control of the form 
\begin{equation}\label{NSC}
	\hat c_i(t)=\frac{\E\left[\E[\zeta_i(t)|\F_t]\middle|\calE_t^{(i)}\right]}{\E\left[\E[p_i(t)|\F_t]\middle|\calE_t^{(i)}\right]}=\frac{\E\left[\zeta_i(t)\middle|\calE_t^{(i)}\right]}{\E\left[p_i(t)\middle|\calE_t^{(i)}\right]},
\end{equation}
is optimal for the maximization problem \eqref{optimalRecursiveF} under $\bbF$ as we have that $\partial_{c_i} \calH^{\bbF}_i \big|_{c_i=\hat c_i}=0$. At this point, if condition \eqref{condition p} is satisfied, we are able to obtain an explicit representation for $\hat c_i(t)$ simply by substituting \eqref{soluzione_esplicita_X}, \eqref{solution_z(t)} and \eqref{soluzione_esplicita_p} in \eqref{NSC}.

\subsection{A zero-sum utility consumption game}
Inspired by the work done in \cite{Tusheng-Oksendal-Sulem} for stochastic delayed equations with constant delay, we consider an optimal consumption problem. However we depart from that work by considering the whole past history in the forward dynamics. Moreover, we take a zero-sum game between two players competing against each other, instead of a coupled delay equation.\\
 We assume that $\alpha$ and $\gamma$ satisfy the necessary integrability conditions and take $\alpha$ to be a bounded deterministic function with $\partial_t\alpha(t,s) \in L^2(dt)$ and $\gamma$ to be a bounded adapted process with $\int_\R\gamma^2(t,z)\nu(dz)<\infty$. Consider a cash flow $X$ with dynamics for $\t$ given by
\begin{equation}\label{Xt example1}
	X(t)=X_0+\int_0^tX(s)\alpha(t,s)-c_1(s)-c_2(s)ds -\int_0^t\!\!\!\int_\R (c_1(s)+c_2(s))\gamma(s,z)\mu(dsdz),
\end{equation}
and $X_0\in\R$. We suppose that, at time $\t$, we consume at the rate $c_1(t)+c_2(t)$, for $c_i(t)\geq 0$, $i=1,2$ càdlàg adapted processes representing the consumption of each player. Assume also that a solution to \eqref{Xt example1} exists.
For each $i=1,2$, we also consider the backward dynamics \eqref{Y_t} given by 
\begin{equation}
	Y_i(t)=h_i(X(T))-\int_t^T\eta_i(s)Y_i(s)ds+\int_t^T\Theta_i(s,z)\mu(dsdz)+\int_t^T dM_i(s), \quad \t,
\end{equation}
for some functions $h_i:\R\longrightarrow\R$ and $\eta_i:[0,T]\longrightarrow \R$ satisfying the necessary integrability conditions that ensure us the existence of a solution $(Y_i, \Theta_i, M_i)$, $i=1,2$.
We want to find the optimal consumption rate $\hat{c}=(\hat c_1, \hat c_2)$ in a zero-sum game such that, for each player $i=1,2$,
\begin{equation*}
	J_i(\hat{c})=\sup_{c\in\calA^{\calE}}J_i(c),
\end{equation*} 
 where 
\begin{equation*}
	J_i(c) := J_i(c_1,c_2)=\E\left[\int_0^TF^i(t,c_1,c_2)dt+V_i(X(T)) \right].	
\end{equation*}
For $i=1,2$, $F^i(t,c_1,c_2):[0,T]\times(\R^+)^2\times \Omega\longrightarrow \R$ is a stochastic utility function such that 
$t\longmapsto F^i(t,c_1,c_2)$ is $\bbF$-adapted for each $c$, and $C^1$ with respect to $c_1$, $c_2$ and $V_i:\R\longrightarrow\R$ is a given concave function $C^1$ in $x$. We assume that at time $\t$ both players access the same information $\E^{(1)}=\E^{(2)}=\bbF$.
In the zero-sum stochastic game setting we have that $J_1(c_1, c_2)=-J_2(c_1,c_2)$. Hereafter we study the maximisation of $J(c_1,c_2):=J_1(c_1,c_2)$ according to Section 3.2.\\
The Hamiltonian functional \eqref{HamiltonZeroF} takes the form
\begin{align}
	\calH^{\bbF}(t)&=F(t,c_1(t), c_2(t))+\big(\alpha(t,t)X(t)- (c_1(t)+c_2(t))\big)\E[p(t)|\F_t]\nonumber \\
	&\quad - \gamma(t,0)\E[q(t,0)](c_1(t)+c_2(t))\lambda_t^B
	 - \int_\R\gamma(t,z)\E[q(t,z)|\F_t](c_1(t)+c_2(t))\nu(dz)\lambda_t^H \nonumber \\
	 &\quad + \eta(t)Y(t)\E[\zeta(t)|\F_t] + \int_0^t\partial_t\alpha(t,s)X(s) ds \ \E[p(t)|\F_t]\label{HamiltonExampleF}
\end{align}
 and the adjoint system \eqref{dp_i}-\eqref{dZ_i} in the setting of \eqref{adjoint_zero_sum} can be rewritten as
\begin{align}
	dp(t)&=\Bigg[\alpha(t,t) + \int_0^t\partial_t\alpha(t,s)ds\Bigg]p(t) dt+\int_\R q(t,z)\mu(dtdz), \label{dp}\qquad \t,\\
	p(T)&=-\partial_xV(X(T))+h(X(T))\zeta(t),\nonumber
\end{align}
and
\begin{equation}
	d\zeta(t)=\eta(t)\zeta(t)dt, \qquad \zeta(0)=0 \label{dZ},
\end{equation}
from which we have $\zeta(t)\equiv 0$, which gives $p(T)=-\partial_xV(X(T))$ in \eqref{dp}. From \cite{DiNunno1} Theorem 4.5 we see that

\begin{equation}\label{explicitsolution}
	(p, q) = \left(-\partial_x V(X(T))\exp\left\{\int_t^T{\bar\alpha}(s)ds \right\},0\right),
\end{equation}
where $\bar\alpha(t) = \alpha(t,t)+\int_0^t\partial_t\alpha(t,s)ds$ is the solution to \eqref{dp}. As for the solution to \eqref{Xt example1} we can simply use the same approach we used for \eqref{soluzione_esplicita_X}.  

We claim now that a control $c_1\geq 0$ satisfying
\begin{align*}
	\partial_{c_1} F(t,c_1,c_2)&=\E\left[\partial_xV(X(T))\exp\left\{\int_t^T{\bar\alpha}(r,t)dr \right\} \middle|\F_t\right]\\
	&=\E\Big[\partial_xV(X(T))\Big|\F_t\Big]\exp\left\{\int_t^T{\bar\alpha}(r,t)dr \right\}
\end{align*}
is optimal and similarly for $c_2$. In this case, in fact, the conditions for the sufficient maximum principle are satisfied as, from Theorem \ref{SufficientMaximumZeroSumTheoremF}, we have that $\partial_{c_1} \calH_1^\bbF$ has to be equal to 0 for a control to be optimal. In particular, this means that we have to ask that
\begin{equation*}
	\partial_{c_1} F(t,c_1,c_2)=\E\left[p_1(t)|\F_t\right],
\end{equation*}
which is exactly what we have above. Notice that in this case we did not make use of the Transformation Rule, as $dX(t)$ can simply be written as
\begin{align*}
    dX(t) &= \Bigg(\alpha(t,t)X(t) - (c_1(t)+c_2(t)) +\int_0^t\partial_t\alpha(t,s)X(s)ds\Bigg)dt\\
    &\quad + \int_\R\gamma(t,z)(c_1(t)+c_2(t))\mu(dt,dz).
\end{align*}

\section{Deterministic time-change and BSVIEs}

    Notice that, when the time-change rates $\lambda^B$ and $\lambda^H$ are deterministic, we could have considered a wider class of backward equations for \eqref{Y_t}. In this case, in fact, one has that $\bbF = \bbG$. In this simpler framework, instead of considering a BSDE such as \eqref{Y_t}, we could have a backward stochastic Volterra integral equation (BSVIE) of the form
    \begin{align}
	    Y_i(t)&=h_i(X(T))-\int_t^Tg_i(t,s,\lambda_s,u(s),X(s\shortminus),Y_i(s\shortminus),\Theta_i(t,s,\cdot))ds\nonumber\\
	    &\quad +\int_t^T\!\!\!\int_\R\Theta_i(t,s,z)\mu(dsdz). \quad \t,\label{Y_t Volterra}
    \end{align}
    in the study of Nash equlibria with \ref{J}. The study of solutions of $(Y_i,\Theta_i) \in L^2(dt\times dP)\times L^2(d\Lambda\times dt \times dP)$ has been carried out in the paper \cite{Popier}. Working with $\bbF=\bbG$, we do not have the orthogonal martingale term $M_i$ which appears in \eqref{Y_t}. 
    We remind that, even when $\lambda_s^B$ and $\lambda_s^H$ are deterministic, the time change \eqref{timechange} is not necessarily a Lévy subordinator.
    
    The first step in the direction of proving a sufficient maximum principle such as Theorem \ref{SufficientMaximumTheoremF} would thus be to define new Hamiltonian functionals accommodating for the specifics of the Volterra structure of \eqref{Y_t Volterra}. We then define the functionals:
    \begin{equation*}
	\calH_i:[0,T]\times\Xi_{\R_+^2}\times \Xi_{\R}\times \Xi_{\R} \times \Xi_{\mathcal Z}\times\Xi_{\bbU}\times\R\times\calR^\bbG\times\calR^\bbG\times\mathcal{I}^\bbG\times\Omega\longrightarrow\R,
    \end{equation*}
    where 
    \begin{align*}
        &\calH_i (t,\lambda,x,y_i,\theta_i,u, y_i^0, \zeta_i,p_i,q_i)\\
	&:=H_0^i(t,\lambda,x,y_i,\theta_i,u,y_i^0,\zeta_i,p_i,q_i)+H_1^i(t,\lambda,x,y_i,\theta_i,u,y_i^0,\zeta_i,p_i,q_i),\nonumber
    \end{align*}
    and $H^i_0$ and $H^i_1$ are defined as
    \begin{align}
	&H_0^i(t,\lambda,x,y_i,\theta_i,u,y_i^0,\zeta_i,p_i,q_i)\\
	&\quad:=F_i(t,u,x(t),y_i(t))+b(t,t,\lambda_t,u(t),x(t))p_i+\kappa(t,t,0,\lambda_t,u(t),x(t))q_i(t,0)\lambda_t^B\nonumber\\
	&\qquad+\int_{\R_0}\kappa(t,t,z,\lambda_t,u(t),x(t))q_i(t,z)\lambda_t^H\nu(dz)+g_i(t,t,\lambda,u(t),x(t),y_i(t),\theta(t,s,\cdot))\zeta_i(t)\nonumber\\
	&H_1^i(t,\lambda,x,y_i,\theta_i,u,y_i^0,\zeta_i,p_i,q_i)\\
	&\quad:=\int_0^t\partial_t b(t,s,\lambda_s,u(s),x(s))ds\ p_i(t)+\int_0^t\!\!\!\int_{\R}\partial_t\kappa(t,s,z,\lambda_s,u(s),x(s))\calD_{s,z}p_i(t)\Lambda(dsdz)\nonumber\\
	&\qquad +\Bigg(\int_t^T\partial_t g_i(t,s,\lambda_s,u(s),x(s),y_i(s),\theta(t,s,\cdot))ds\nonumber\\
	&\qquad+\int_t^T\partial_{\theta_0}g_i(t,s,\lambda_s,u(s),x(s),y_i(s),\theta(t,s,\cdot))\partial_t\theta(t,s,\cdot)\lambda_s^B ds\nonumber\\
	&\qquad+\int_t^T\left\langle\nabla_{\theta_{z}}g_i(t,s,\lambda_s,u(s),x(s),y_i(s),\theta(t,s,\cdot)),\partial_t\theta(t,s,\cdot)\right\rangle \lambda_s^H ds \Bigg) \zeta_i(t),\label{Hamilton_Volterra}
\end{align}
    where $(p_i, q_i)$ solve \eqref{dp_i} and $\zeta_i$ is a solution to \eqref{dZ_i} and where we denoted the action of the operator $f$ on the function $g$ by $\langle f, g\rangle$.
    
    For later application of a Transformation Rule similar to the one presented in \cite{Protter} for forward Volterra integral equations, we need some regularity conditions. Hence we assume that, for all $z\in\R$ the partial derivative of $\Theta$ with respect to $t$ is locally bounded (uniformly in $t$) and satisfies
	\begin{equation}\label{trasformation rule hp back}
	    |\partial_t\Theta(t_1,s,z)- \partial_t \Theta(t_2,s,z)|\leq K |t_1-t_2|,
	\end{equation}
	for some $K>0$ and for each fixed $s\leq t_1, t_2$, $z\in\R$.
	In this case we have the following:
		\begin{lema}(Transformation Rule for backward Volterra integral equations)\label{Transformation rule back}
The backward equation \eqref{Y_t Volterra} satisfying condition \eqref{trasformation rule hp back}, can be rewritten in differential notation as 
\begin{align*}
dY_i(t) &= - \Bigg(g_i(t,t,\lambda_t,u(t),X(t), Y_i(t),\Theta_i(t,t,\cdot)) \\&\quad+\int_t^T\partial_t g_i(t,s,\lambda_s,u(s),X(s\shortminus),Y_i(s\shortminus),\Theta_i(t,s,\cdot))ds\\
	    & \quad-\int_t^T\!\!\!\int_\R\partial_t\Theta_i(t,s,z)\mu(dsdz)\Bigg)dt + \int_\R\Theta_i(t,t,z)\mu(dtdz).
\end{align*}
\end{lema}
\begin{proof}
    The proof of this statements follows the one presented in \cite{Protter} for forward Volterra integral equations and in \cite{DiNunnoGiordano} in the context of time change. We report it here for completeness. We have that, using the shortened notation
    $$g(t,s):= g(t,s,\lambda_s,u(s), X(s), Y(s), \Theta(t,s,\cdot))$$
    and dropping the superscript $i$,
    \begin{align*}
        Y(t) &= -\int_t^T g(t,s)ds +\int_t^T\!\!\!\int_\R\Theta(t,s,z)\mu(dsdz)\\
        &=-\int_t^T g(t,s)ds -\int_t^T\!\!\!\int_\R\Theta(t,s,z)-\Theta(s,s,z)\mu(dsdz)+\int_t^T\!\!\!\int_\R\Theta(s,s,z)\mu(dsdz)
    \end{align*}
    Noticing now that $\Theta(t,s,z)-\Theta(s,s,z) = -\int_t^s\partial_r\Theta(r,s,z)dr = -\int_t^T\mathds{1}_{[t,s]}(r)\partial_r \Theta(r,s,z) dr$, $s>t$ we can apply the Fubini theorem for stochastic integration as in \cite{Jacod} and we obtain that
    \begin{align*}
        \int_t^T\!\!\!\int_\R \Theta(t,s,z)-\Theta(s,s,z)\mu(dsdz) &= -\int_t^T\!\!\!\int_\R\left\{\int_t^T\mathds{1}_{[t,s]}(r)\partial_r \Theta(r,s,z) dr\right\}\mu(dsdz)\\
        &=-\int_t^T\left\{\int_r^T\!\!\!\int_\R\partial_r\Theta(r,s,z)\mu(dsdz)\right\}dr.
    \end{align*}
    The good definition and Lebsesgue integrability of $\int_r^T\!\!\!\int_\R\partial_r\Theta(r,s,z)\mu(dsdz)$ comes now from \cite{Protter} Theorem 3.2. From here we conclude. 
\end{proof}

    In the setting of the current section, a sufficient maximum principle along the lines of Theorem \ref{SufficientMaximumTheoremF}, but with Volterra dynamics for $Y$ and Hamiltonian $\calH^{\bbF}$ defined as $\calH^{\bbF}:=\E[\calH(t)|\F_t]$ for $\calH$ in \eqref{Hamilton_Volterra} holds. The statement would read as follows:
    \begin{teo}\label{SufficientMaximumTheoremF_Volterra}
	Let $\hat u= (\hatu_1,\hatu_2)\in\calA^{\calE}_1\times\calA^{\calE}_2 = \calA^{\caltildeE}_1\times\calA^{\caltildeE}_2$ and assume that the corresponding solutions $\hatX$, $(\hatY_i, \hat\Theta_i, \hat M_i)$, $(\hatp_i, \hatq_i)$, $\hat\zeta_i$ of equations \eqref{X_t}, \eqref{Y_t Volterra}, \eqref{dp_i}, and \eqref{dZ_i} (with $y_i^0 = \hatY_i(0)$ in the initial value) exist for $i=1,2$. 
	We assume that the functionals $h_i$ in \eqref{Y_t Volterra} are concave. We also consider performance functionals \eqref{J} with concave $\varphi_i$ and $\psi_i$.
	Assume that, for all $i=1,2$, $\calH^\bbF_i$ in \eqref{HamiltonF} satisfy the following conditions:
	\begin{enumerate}
	\item[(i)]
	The maps
		$$x,y,\theta \longmapsto \calH_i^\bbF(t,\lambda,x,y,\theta,u,y^0,\zeta,p,q)$$
		are concave for all $\t$, $\lambda \in \calL$, $u\in\bbU$, $y^0\in \R$, $z,p \in\calR^\bbG$, $q\in\mathcal I^\bbG$.\\
	\item[(ii)] The following extremes are achieved:
		\begin{align*}
		\sup_{v\in\calA^{\calE}_1}&\E\left[\calH^{\bbF}_1(t,\lambda,\hatX,\hatY_1,\hat\Theta_1(t,\cdot),v,\hatu_2,\hat \zeta_1,\hatp_1,\hatq_1)|\calE_t^{(1)}\right]\\
		&=\E\left[\calH^{\bbF}_1(t,\lambda,\hatX,\hatY_1,\hat\Theta_1(t,\cdot),\hatu_1,\hatu_2,\hat \zeta_1,\hatp_1,\hatq_1)|\calE_t^{(1)}\right],
		\end{align*}
		and
		\begin{align*}
		\sup_{v\in\calA^{\calE}_2}&\E\left[\calH^{\bbF}_2(t,\lambda,\hatX,\hatY_2,\hat\Theta_2(t,\cdot),\hatu_1,v,\hat \zeta_2,\hatp_2,\hatq_2)|\calE_t^{(2)}\right]\\
		&=\E\left[\calH^{\bbF}_2(t,\lambda,\hatX,\hatY_2,\hat\Theta_2(t,\cdot),\hatu_1,\hatu_2,\hat \zeta_2,\hatp_2,\hatq_2)|\calE_t^{(2)}\right].
		\end{align*}
		\end{enumerate}
		Then $\hatu=(\hatu_1,\hatu_2)$ is a Nash equilibrium with respect to the information flow $\bbF=\bbG$ for the stochastic game \eqref{X_t},\eqref{Y_t},\eqref{J}.
\end{teo}
    
    The proof relies on writing $Y$ in differential notation by applying the Transformation Rule above and then using the Itô formula for the product in \eqref{I3_inequality} to compute $\E[\hat\zeta_1(T)(Y_1(T)-\hat Y_1(T))]$.

	It is not obvious to what extent condition \eqref{trasformation rule hp back} hold for the solution of general BSVIEs. Inspired by \cite{OksendalForward}, where the authors prove that such smoothness conditions hold for a linear BSVIE with Lévy noise, we consider the case when when $g_i(t,s,\lambda_s,u(s),X(s),Y_i(s),\Theta_i(t,s,\cdot))$ is linear in $Y_i$ and study explicit solutions for \eqref{Y_t Volterra}. 
    
\subsection{Explicit solution of a class of linear BSVIE}
	
    In this subsection we consider the BSVIE of the following type
	\begin{equation}\label{Y_t_linear}
		Y(t)=\xi(t)+\int_t^T \gamma(t,s,\lambda_s)Y(s) ds -\int_t^T\!\!\!\int_\R\Theta(t,s,z)\mu(dsdz), \quad \t,
	\end{equation}
	where the driver of \eqref{Y_t} is reduced to $$g(t,s,\lambda(s),Y(s),\Theta(t,s,\cdot))=\gamma(t,s,\lambda(s))Y(s),$$ with $\gamma(t,s,\lambda(s))$, $s\in[t,T]$, bounded $\bbF^\Lambda$-adapted process,  for all $\t$ and $\xi(t):[0,T]\times \Omega\longrightarrow\R$.	
	For any $\t$, we define recursively the process:
	\begin{align*}
		\gamma^{(1)}(t,r,\lambda_r)&:=\gamma(t,r,\lambda_r), \quad r\in[t,T],\\
		\gamma^{(n)}(t,r,\lambda_r)&:=\int_t^r\gamma^{(n-1)}(t,s,\lambda_s)\gamma(s,r,\lambda_r)ds, \quad r\in[t,T], \quad n\geq 2
	\end{align*}
	and we put \begin{equation}\label{defPsi}
	\Psi(t,r,\lambda_r):=\sum_{n=1}^\infty \gamma^{(n)}(t,r,\lambda_r), \ 0\leq t \leq r\leq T.
	\end{equation}
	The series converges in $L^2(dP)$ since $|\gamma(t,r,\lambda_r)|\leq C$ for some constant $C>0$. In fact, by induction we can see that this condition implies that $|\gamma^{(n)}(t,r,\lambda_r)|\leq\frac{C^nT^n}{n!}$, for all $t,r,n$. 
	\begin{teo}\label{ExplicitY}
		Let $\Psi$ in \eqref{defPsi} be well defined in $L^2(dP)$. The linear BSVIE \eqref{Y_t_linear} admits an explicit solution $(Y,\Theta)\in L^2(dt\times dP)\times \mathcal{I}^{\bbG}$, where
		\begin{itemize}
		\item[(i)] the component $Y$ is given by
		\begin{equation}\label{explicit_Y}
			Y(t)=\E[\xi(t)|\calG_t]+\int_t^T\Psi(t,r,\lambda_r)\E[\xi(r)|\calG_t]dr,
		\end{equation}
		\item[(ii)] the component $\Theta$ is given by the NA-derivative
		\begin{equation}\label{explicit_phi}
			\Theta(t,s,z)=\calD_{s,z}U(t), \ 0\leq t \leq s \leq T, \ z\in\R,
		\end{equation}
		of
		\begin{equation}\label{defU}
			U(t):=\xi(t)+\int_t^T \gamma(t,r,\lambda_r)Y(r)dr-Y(t), \quad 0\leq t\leq T.
		\end{equation}
		\end{itemize}
	\end{teo}
\begin{proof}
	For $u\leq t$, we take the conditional expectation given $\calG_u$ of $Y(t)$ in \eqref{Y_t_linear} 
	\begin{equation}\label{OK3_19}
		\E[Y(t)|\calG_u]=\E[\xi(t)|\calG_u]+\int_t^T\gamma(t,s,\lambda_s)\E[Y(s)|\calG_u]ds.
	\end{equation}	
	Denote $\tilde{Y}(t,u)=\E[Y(t)|\calG_u]$ for $u\leq t$, and $\tilde{\xi}(t,u):=\E[\xi(t)|\calG_u]$, then \eqref{OK3_19} can be rewritten as
	\begin{equation*}
		\tilde{Y}(t,u)=\tilde{\xi}(t,u)+\int_t^T\gamma(t,s,\lambda_s)\tilde{Y}(s,u)ds, \quad u\leq t\leq T.
	\end{equation*}
	By substituting $\tilde{Y}(s,u)=\tilde{\xi}(s,u)+\int_s^T\gamma(s,r,\lambda_r)\tilde{Y}(r,u)dr$ in the previous equation, we obtain that
	\begin{align*}
		\tilde{Y}(t,u)&=\tilde{\xi}(t,u)+\int_t^T \gamma(t,s,\lambda_s)\left(\tilde{\xi}(s,u)+\int_s^T\gamma(s,r,\lambda_r)\tilde{Y}(r,u)dr\right)ds\\
		&=\tilde{\xi}(t,u)+\int_t^T\gamma(t,s,\lambda_s)\tilde{\xi}(s,u)ds+\int_t^T\gamma^{(2)}(t,r,\lambda_r)\tilde{Y}(r,u)dr,
	\end{align*}
	by Fubini theorem. Hence, by iteration, we obtain
	\begin{align}
		\tilde{Y}(t,u)&=\tilde{\xi}(t,u)+\int_t^T\sum_{n=1}^N\gamma^{(n)}(t,r,\lambda_r)\tilde{\xi}(r,u)dr+\int_t^T\gamma^{(N+1)}(t,r,\lambda_r)\tilde{Y}(r,u)dr.\label{NUM back sol}	
	\end{align}
Thanks to the Hölder inequality, the bound on $\gamma^{(N+1)}$, and the Jensen inequality, we have
\begin{align*}
	\E\left[\left(\int_t^T\gamma^{(N+1)}(t,r,\lambda_r)\tilde{Y}(r,u)dr\right)^2\right]&\leq \E\left[\int_t^T\gamma^{(N+1)}(t,r,\lambda_r)^2dr\int_t^T\tilde{Y}(r,u)^2dr\right]\\
	&\leq T\left(\frac{C^{N+1}T^{N+1}}{(N+1)!}\right)\E\left[\int_t^T\E\left[Y(r)^2|\calG_u\right]dr\right]\\
	&\leq T\left(\frac{C^{N+1}T^{N+1}}{(N+1)!}\right)\int_t^T\E[Y(r)^2]dr\longrightarrow 0,
\end{align*}
	as $N$ goes to infinity. Hence taking the limit in $L^2(dP)$ in \eqref{NUM back sol}, we obtain
	\begin{align*}
		\tilde{Y}(t,u)&=\tilde{\xi}(t,u)+\int_t^T\sum_{n=1}^\infty\gamma^{(n)}(t,r,\lambda_r)\tilde{\xi}(r,u)dr\\
		&=\tilde{\xi}(t,u)+\int_t^T\Psi(t,r,\lambda_r)\tilde{Y}\tilde{\xi}(r,u)dr
	\end{align*}
	where $\Psi$ defined in \eqref{defPsi}.
 Let us now take $u=t$ in the equation above and \eqref{OK3_19}, then we have 	 
	 $$Y(t)=\E[Y(t)|\calG_t]=\tilde{Y}(t,t),$$
	 which yields $(i)$. 	
	Given $(i)$, $U(t)$, $\t$ is well defined. Also observe that \eqref{Y_t_linear} yields
	\begin{equation*}
		U(t)=\int_t^T\!\!\!\int_\R\Theta(t,s,z)\mu(dsdz).
	\end{equation*}
	Hence, by Theorem \ref{C-Ok-nonanticipating} we have that, for all $\t$, $\Theta(t,s,z)=\calD_{s,z}U(t)$, $s\in[t,T]$, $z\in\R$.	
\end{proof}

We can now go back to look for conditions on $\Theta(t,s,z) = \calD_{s,z}U(t)$ that ensure us that $\partial_t\Theta(t,s,z)$ exists and is Lipschitz.
By substituting \eqref{explicit_phi} and \eqref{explicit_Y} in \eqref{defU} we find that that:
\begin{align*}
    \partial_t\Theta(t,s,z) &= \partial_t\calD_{sz}\Bigg(\xi(t) + \int_t^T\gamma(t,r,\lambda_r) Y(r)dr - Y(t)\Bigg)\\
    &= \partial_t\calD_{sz}\Bigg(\xi(t) +\int_t^T\gamma(t,r,\lambda_r)\left\{\E[\xi(r)|\calG_r]+\int_r^T\Psi(r,s,\lambda_s)\E[\xi(s)|\calG_s]ds\right\}dr\\ 
    &-\E[\xi(t)|\calG_t]+\int_t^T\Psi(t,s,\lambda_s)\E[\xi(s)|\calG_s]ds\Bigg).
\end{align*}
    From here we see that, if $\xi$ and $\gamma$ are $C^1$ with respect to $\t$ (uniformly for $s\in[0,T]$, $\lambda\in[0,\infty)^2$) with derivatives Lipschitz with respect to $t$, then the required condition on $\Theta$ to apply thr Transformation Rule hold.

\subsection{A linear utility consumption example}
        Having ensured that the required conditions to apply the Transformation Rule hold, we can present an example similar to the one in subsection 4.1 but with Volterra dynamics for the backward \eqref{Y_t}. We take once again  $X(t)$ to have dynamics
    \begin{equation}\label{X_t Volterra}
        X(t) = X_0 +\int_0^t \alpha(t,s)X(s) - c_1(s)-c_2(s)ds - \int_0^t\!\!\!\int_\R\gamma(s,z) (c_1(s)+c_2(s)) \mu(dsdz)
    \end{equation}
    and $Y_i(t)$ follows the dynamics
\begin{equation}
	Y_i(t)=X(T)+\int_t^T\eta_i(t,s,c_i(s))Y(s)ds-\int_t^T\Theta_i(t,s,z)\mu(dsdz), \quad \t,
\end{equation}
for $X_0\in\R$ and some functions $\alpha:[0,T]^2\longrightarrow\R$, $\gamma:[0,T]\times\R\times\Omega\longrightarrow\R$ and $\eta:[0,T]^2\times \bbU\times\Omega\longrightarrow \R$ satisfying the necessary integrability conditions, and such that $\alpha$ is a bounded deterministic function with $\partial_t\alpha(t,s) \in L^2(dt)$ and $\gamma$ is a bounded adapted process with $\int_\R\gamma^2(t,z)\nu(dz)<\infty$. As we saw in the previous section, we assume $\eta_i$ to be $C^2$ with respect to $\t$ and such that $\partial_t\eta_i$ is Lipschitz continuous with respect to $t$ for all $i=1,2$. We also assume that $\eta_i(t,s,c_i(s)) $ is bounded with derivative $\partial_t\eta_i(t,s,c_i(s))\in L^2(dt\times dP)$ for $i=1,2$. From the examples in the previous section we already know that in this case we do not need a Transformation Rule for $X$. Once again we suppose we are in the zero-sum case and we will drop the superscript $i$ whenever no confusion arises.
We want to find the optimal consumption rate $\hat{c}=(\hat c_1, \hat c_2)$ such that
\begin{equation*}
	J(\hat{c})=\sup_{c\in\calA^{\bbF}}J(c),
\end{equation*} 
 where 
\begin{equation*}
	J(c):=J(c_1,c_1)=\E\left[\int_0^TF(t,c_1(t),c_2(t),\omega)dt+K(( X(T)+Y(0)) \right],
\end{equation*}
$F$ is $C^1$ with respect to $c$ and $K \in \R_0$.
We start by rewriting the Hamiltonian functional \eqref{Hamilton_Volterra}:
\begin{align}
	\calH_i(t)&=F(t,c_1,c_2)+[\alpha(t,t)X(t)-c_1(t)-c_2(t)]p_i(t)+ \int_0^t \partial_t\alpha(t,s)X(s)ds\ p_i(t)\nonumber \\
	&\quad -\gamma(t,0)q_i(t,0)\lambda^B_t(c_1(t)+c_2(t))-\int_{\R_0}\gamma(t,z)q_i(t,z)(c_1(t)+c_2(t))\lambda^H_t\nu(dz)\nonumber\\
	&\quad +\Bigg(\eta_i(t,t,c_i(t))Y_i(t) + \int_t^T \partial_t\eta_i(t,s,c_i(s))Y_i(s) ds\Bigg) \zeta_i (t) 
\end{align}
In this case we have that the forward adjoint equation \eqref{dZ_i} can be written as 
\begin{equation}\label{Z_t_Volterra}
    d\zeta_i(t) = \Bigg(\eta(t,t,c_i(t))+\int_t^T\partial_t\eta(t,s,c_i(s)) ds \Bigg) \zeta_i(t)dt, \quad \zeta_i(0)=K
\end{equation}
whereas the backward adjoint equation \eqref{dp_i} is \begin{align}
	dp(t)&=\Bigg[\alpha(t,t)+\int_0^t\partial_t\alpha(t,s)ds\Bigg]p(t)dt
	\int_\R q(t,z)\mu(dtdz),\nonumber \\
	p(T)&=K\label{dp2_v}.
\end{align}
Analogously to \eqref{explicitsolution} we have that the solution for \eqref{dp2_v} is given by 
\begin{equation}\label{soluzione esplicita p esempio 3}
   (p,q) =\Bigg( K\exp\left\{\int_0^t\bar\alpha(s)ds\right\},0\Bigg)
\end{equation}
where $\bar\alpha(t) = \alpha(t,t)+\int_0^t\partial_t\alpha(t,s)ds$.\\
Similarly for what we did in \eqref{soluzione_esplicita_X} we find that the explicit solution for \eqref{X_t Volterra} is given by 
\begin{equation}\label{soluzione esplicita X esempio 3}
    X(t) = M(t) \left(X_0-\int_0^tc_1(s)+c_2(s) ds - \int_0^t\!\!\!\int_\R\gamma(s,z)(c_1(s)+c_2(s)) \mu(dsdz)\right),
\end{equation}
where $M(t)$ is defined in \eqref{M(t)}.

 As for the linear BSVIE, using Theorem \ref{explicit_Y} we have that a solution is given by 
 \begin{align}
     Y(t) &= \E[X(T)|\calG_t]+\int_t^T\Psi(t,r,c)\E[X(T)|\calG_t],\nonumber\\
     \Theta(t,s,z) &=\calD_{s,z}U(t),\label{soluzione esplicita Y esempio 3}
 \end{align}
 where $\Psi$ is defined in \eqref{defPsi} and $U(t)$ in \eqref{defU}. 
 Also, one can find that the solution to  \eqref{Z_t_Volterra} is given by
 \begin{equation}\label{soluzione esplicita Z esempio 3}
     \zeta_i(t) = K\exp\left\{-\int_0^t\bar\eta(s,c_i(s))ds\right\}
 \end{equation}
 where $\bar\eta_i(t,c_i(t)) = \eta(t,t,c_i(t))+\int_t^T\partial_t\eta(t,s,c_i(s))ds$. 
 We now claim that a positive control $c_i$ that satisfies
 \begin{align}
     0&= \partial_{c_i} F(t,c_1,c_2) - p_i(t)\nonumber \\
     &\quad -\gamma(t,0)q_i(t,0)\lambda_t^B - \int_{\R_0}\gamma(t,z)q_i(t,z)\lambda_t^H\nu(dz)\nonumber\\
     &\quad+\partial_{c_i} \Bigg(\eta(t,t,c_i(t))Y(t) + \int_t^T \partial_t\eta(t,s,c_i(s))Y_i(s) ds\Bigg) \zeta_i(t) \label{SuffVolterra}
 \end{align}
 is optimal, $i=1,2$. Having solved both the forward-backward Volterra integral equations system and the adjoint equations systems, we can substitute \eqref{soluzione esplicita p esempio 3}, \eqref{soluzione esplicita X esempio 3}, \eqref{soluzione esplicita Y esempio 3}, \eqref{soluzione esplicita Z esempio 3}, them into \eqref{SuffVolterra} and obtain our sufficient condition for optimality.

	\subsection*{Acknowledgment} The research leading to these results is within the project STORM: Stochastics for Time-Space Risk Models, number: 274410, of the Research Council of Norway.

\bibliographystyle{plain}
\bibliography{bib}
\end{document}